\documentclass[aos,MSNbibl,nameyear,dvips]{arximspdf}
\usepackage{dcolumn}
\usepackage{graphicx}

\doi{10.1214/11-AOS925}
\volume{39}
\issue{6}
\pubyear{2011}
\firstpage{2912}
\lastpage{2935}

\makeatletter

\newcolumntype{d}[1]{D{.}{.}{#1}}
\newcolumntype{k}[1]{D{,}{}{#1}}

\newtheorem{theorem}{Theorem}[section]
\newtheorem{lemma}{Lemma}[section]

\newcommand{\bbf}{\mathbf{f}}
\newcommand{\by}{\mathbf{y}}
\newcommand{\bz}{\mathbf{z}}
\newcommand{\var}{\operatorname{var}}

\newcommand{\bff}{\mathbf{f}}

\newcommand{\btheta}{{\bolds\theta}}
\newcommand{\eps}{\varepsilon}
\newcommand{\bmu}{{\bolds{\mu}}}
\newcommand{\bR}{\mathbf{R}}
\newcommand{\bAIC}{\operatorname{AIC}}
\newcommand{\bBIC}{\operatorname{BIC}}
\newcommand{\bMDL}{\operatorname{MDL}}

\newcommand{\bx}{\mbox{$x$}}

\newcommand{\argmin}{\arg\min}

\makeatother

\begin{document}
\begin{frontmatter}

\title{On image segmentation using information theoretic~criteria}
\runtitle{Consistent image segmentation}

\begin{aug}
\author[A]{\fnms{Alexander} \snm{Aue}\thanksref{t1}\ead[label=e1]{alexaue@wald.ucdavis.edu}}
\and
\author[A]{\fnms{Thomas C. M.} \snm{Lee}\corref{}\thanksref{t2}\ead[label=e2]{tcmlee@ucdavis.edu}}
\runauthor{A. Aue and T. C. M. Lee}
\affiliation{University of California at Davis}
\address[A]{Department of Statistics \\
University of California at Davis \\
4118 Mathematical Sciences Building \\
One Shields Avenue \\
Davis, California 95616 \\
USA\\
\printead{e1}\\
\phantom{E-mail: }\printead*{e2}} 
\end{aug}

\thankstext{t1}{Supported in part by NSF Grant 0905400.}
\thankstext{t2}{Supported in part by NSF Grants 0707037 and 1007520.}

\received{\smonth{6} \syear{2011}}
\revised{\smonth{9} \syear{2011}}

%
\begin{abstract}
Image segmentation is a long-studied and important problem in image
processing. Different solutions have been proposed, many of which
follow the information theoretic paradigm. While these information
theoretic segmentation methods often produce excellent empirical
results, their theoretical properties are still largely unknown. The
main goal of this paper is to conduct a rigorous theoretical study
into the statistical consistency properties of such
methods. To be more specific, this paper investigates if these
methods can accurately recover the true number of segments together
with their true boundaries in the image as the number of pixels tends
to infinity. Our theoretical results show that both the Bayesian
information criterion (BIC) and the minimum description length (MDL)
principle can be applied to derive statistically consistent
segmentation methods, while the same is not true for the Akaike
information criterion (AIC). Numerical experiments were conducted to
illustrate and support our theoretical findings.
\end{abstract}

%

\end{frontmatter}

\section{Introduction}\label{sec1}
Image segmentation aims to partition an image into a~set of
nonoverlapping regions so that pixels within the same region are
homogeneous with respect to some characteristic (e.g., gray value or
roughness), while pixels from adjacent regions are significantly
different with respect to the same characteristic. It is a
fundamental problem in image processing, as very often it is necessary
to first group the highly localized pixels into more global and
meaningful segmented objects to facilitate the extraction of useful
information. In this paper, gray value is the image characteristic that
forms the basis for segmentation. For general introductions to image
segmentation, see, for example, \citet{Glasbey-Horgan95} and
\citet{Haralick-Shapiro92}.

A grayscale image can be seen as a two-dimensional (2D) surface living
in a three-dimensional space. Therefore one popular approach to
segmenting it is to model it by a 2D piecewise constant function,
with the set of all discontinuity points defining the region
boundaries of the image. Examples of segmentation methods that
follow this approach include
\citet{Kanungo-et-al95},
\citet{LaValle-Hutchinson95}, \citet{Leclerc89},
Lee (\citeyear{Lee98segcor,Lee00segind}), \citet{Luo-Khoshgoftaar06}
and \citet{Wang-et-al09}.
As to be demonstrated below, segmenting images with this approach can
be recast as a model selection problem, and one crucial issue to its
success is the choice of the model complexity, which is equaivalent to
choosing the number of regions together with the shapes of their
boundaries.
Common information theoretic methods such as the Akaike information
criterion (AIC) [\citet{Akaike74}], the Bayesian information criterion
(BIC), also known as the Schwarz information criterion
[\citet{Schwarz78}] and the minimum description length (MDL) principle
[Rissanen (\citeyear{Rissanen89,Rissanen07})] have
been adopted to solve this problem;
for example, see
\citet{Kanungo-et-al95}, \citet{Leclerc89}, Lee
(\citeyear{Lee98segcor,Lee00segind}),
\citet{Luo-Khoshgoftaar06}, \citet{Murtagh-et-al05},
\citet{Stanford-Raftery02}, \citet{Zhang-Modestino90}
and \citet{Zhu-Yuille96}.
While many of these methods produce excellent practical results, their
theoretical properties are still largely unknown. The goal of this
paper is to conduct a systematic study on the theoretical properties
of these methods, with the hope of enhancing our understanding of their
performances, at both theoretical and empirical levels. To the best
of our knowledge, this is the first time that such a rigorous
theoretical study is being performed for image segmentation methods.

The rest of this paper is organized as follows. Background material
is presented in Section \ref{sec2}. Section \ref{secmain} presents our
main theoretical results. These theoretical results are empirically
verified by numerical experiments in Section \ref{sec4}.
Concluding remarks are offered in Section \ref{sec5}, while technical
details are delayed to the \hyperref[app]{Appendix}.

\section{Background}\label{sec2}

Denote by $f$ the true image and $\Xi_n=\{x_1, \ldots, x_n\}$ the set
of $n$ grid points at which a noisy version of $f$ is sampled.
Without loss of generality it is assumed that the domain of $f$ is
$[0,1]^2$. As mentioned before,~$f$ is modeled as a 2D piecewise
constant function as follows. Write $f_i=f(x_i)$ and
$\bbf=(f_1,\ldots,f_n)'$. Let the number of regions (or pieces or
segments) in~$f$ be~$m$, and denote the gray value and domain of the
$\nu$th region as~$\mu_\nu$ and~$R_\nu$, respectively. Then we have,
for $i=1, \ldots, n$,
%
\begin{eqnarray}
\label{eq1a}
f_i &=& \mu_\nu \qquad\mbox{if } x_i \in R_\nu,
\\
%
\label{eq1b}
\bigcup_{\nu=1}^mR_\nu&=&[0,1]^2\quad \mbox{and}\quad
R_\nu\cap R_{\nu^\prime}=\varnothing\qquad\mbox{if } \nu\not=\nu
^\prime.
\end{eqnarray}
In the sequel we write $\bR=(R_1,\ldots,R_m)$ and
$\bmu=(\mu_1,\ldots,\mu_m)'$. Thus $\bR$ defines a segmentation of
$f$. The observed noisy version $\by=(y_1,\ldots,y_n)'$ of~$\bbf$ is
modeled as
%
\begin{equation}\label{eqyi}
y_i=f_i+\eps_i, \qquad  i=1,\ldots,n,
\end{equation}
where the noise $\eps_i$'s are independent, identically distributed
random variables with zero mean and variance $\sigma^2$. Given $\by$,
the goal is then to estimate~$\bbf$, which is equivalent to estimating
$m$, $\bR$ and $\bmu$.

For simplicity, denote by $\btheta_m=(m,\bR,\bmu)'$ a generic
parameter vector. Estimating $\bff$ is hence equivalent to the model
selection problem in which each model is determined by the parameter
$\btheta_m$. Let $\mathrm{RSS}_m=\sum_i(y_i-\hat f_i)^2$ be the
corresponding residual sum of squares. Notice that different values
of $m$ would lead to a different number of parameters in $\btheta_m$.
Also notice that $\btheta_m$ cannot be estimated by minimizing
$\mathrm{RSS}_m$, as $\mathrm{RSS}_m$ can be made arbitrarily small as
$m$ tends to $n$. One way to resolve this issue is to add a penalty
term to $\mathrm{RSS}_m$ to suitably penalize the complexity of
$\btheta_m$. As alluded to before, information theoretic model selection
methods like AIC, BIC and MDL can be used to derive such a penalty.
We first focus on the MDL criterion derived by \citet{Lee00segind},
%
\begin{equation}
\label{eq2}
\bMDL(m,\bR)=m\ln n+\frac{\ln3}{2}\sum_{\nu=1}^{m}b_\nu+
\frac12\sum_{\nu=1}^m\ln a_\nu+\frac n2\ln\biggl(\frac{
\mathrm{RSS}_m}{n}\biggr),
\end{equation}
where each region $R_\nu$ enters through its ``area'' $a_\nu$ (in
terms of number of pixels) and ``perimeter'' $b_\nu$ (in terms of
number of pixel edges). These quantities are formally defined as
\[
a_\nu=\# (\Xi_n\cap R_\nu)  \quad\mbox{and}\quad
b_\nu=\# (\Xi_n\cap\partial R_\nu)
\]
with $\#A$ and $\partial A$ indicating, respectively, cardinality and
boundary of the set $A$. Observe that, once the estimates $\hat{m}$ and
$\hat{\bR}$ are specified, $\bmu$ can be uniquely estimated by
%
\begin{equation}\label{eqmuhat}
\hat{\mu}_\nu=\frac{1}{\hat{a}_\nu}\sum_{i\in\hat{R}_\nu} y_i
\qquad\mbox{for all $\nu$},
\end{equation}
and therefore $\bmu$ is dropped in the argument list of
$\bMDL(m,\bR)$. To sum up, the MDL-based method of
\citet{Lee00segind} estimates $m$ and $\bR$ as the joint minimizer
of~(\ref{eq2}), which is equivalent to saying
%
\begin{equation}
\label{eq3}
(\hat{m},\hat{\bR})=\mathop{\argmin}_{m\leq M,\bR}\frac2n\bMDL(m,\bR),
\end{equation}
and $\hat{\bmu}$ is given by (\ref{eqmuhat}). Practical algorithms,
developed, for example, by \citet{Lee00segind} and
\citet{Zhu-Yuille96},
can be used to solve (\ref{eq3}).

One can also use AIC and BIC to derive penalty terms to add to
$\mbox{RSS}_m$, and the resulting penalties will be\vspace*{1pt} proportional to
the number of ``free'' (and independent) parameters in the fitted\vadjust{\goodbreak}
image $\hat{\bbf}$ [e.g., \citet{Murtagh-et-al05},
\citet{Stanford-Raftery02}
and \citet{Zhang-Modestino90}]. This leads to the following question:
what would
be a meaningful way of counting the number of free parameters in
$\hat{\bbf}$? There seems to be no unique answer, but we shall follow
\citet{Murtagh-et-al05} and \citet{Stanford-Raftery02} and model each
true pixel value $f_i$ with a~mixture distribution of $m$ Gaussians,
where the mean, variance and mixing probability for the $\nu$th
Gaussian are $\mu_\nu$, $\sigma^2$ and $a_\nu/\sum_\nu a_\nu$,
respectively. As there are $m$ of the $\mu_\nu$'s, one $\sigma^2$ and
$m-1$ free mixing probabilities, the total number of free parameters
is $2m$. With this, the corresponding AIC and BIC segmentation
criteria are
\[
\bAIC(m,\bR)=2m+\frac n2\ln\biggl(\frac{\mathrm{RSS}_m}{n}\biggr)
\]
and
\[
\bBIC(m,\bR)=m\ln n +\frac n2\ln\biggl(\frac{\mathrm{RSS}_m}{n}\biggr),
\]
respectively. The AIC and BIC estimates for $(m, \bR)$ are then given
by
%
\begin{equation}\label{eqaicestimate}
(\hat{m},\hat{\bR})=\mathop{\argmin}_{m\leq M,\bR}\frac2n\bAIC(m,\bR)
\end{equation}
and
%
\begin{equation}\label{eqbicestimate}
(\hat{m},\hat{\bR})=\mathop{\argmin}_{m\leq M,\bR}\frac2n\bBIC(m,\bR),
\end{equation}
respectively.
Observe that for both $\bAIC(m,\bR)$ and $\bBIC(m,\bR)$, the region
boundaries $\bR$ are not explicitly penalized; they enter the criteria
only through $\mbox{RSS}_m$. Also observe that the penalty term of
$\bAIC(m,\bR)$ is independent of~$n$.

Before we proceed further, it is worthwhile to point out a major
difference between the variable selection problem in linear regression
models and the image segmentation problem. In variable selection for
linear regression, the goal is to select the significant predictors
and remove the insignificant ones from the model. In other words,
some ``data'' are not used in estimating the model parameters. For
image segmentation, the goal is to group homogeneous pixels together
to form segmented objects, and in this process all data (i.e., all
pixel values) are always used to estimate the model parameters. Given
this major difference, one can see that variable selection in linear
regression and image segmentation are two different problems, and
hence existing theories from classical linear regression modeling
cannot be directly applied to image segmentation.

\section{Main results}
\label{secmain}
This section presents our main theoretical findings. Briefly, both
the BIC and MDL segmentation solutions are statistically consistent in
a well-defined sense, while the AIC solution is not.\vadjust{\goodbreak}

The consistency of the BIC and MDL solutions are investigated at two
levels. First, we will establish the strong consistency of
$\hat{\bR}$ if the true number of regions $m=m^0$ can be assumed
known. Second, if the true value $m^0$ is unknown and if the noise is
restricted to be Gaussian, we will establish the weak consistency of
$\hat{m}$ and $\hat{\bR}$. While the existence of a true underlying
model was not essential for the practical use
of (\ref{eq3})--(\ref{eqbicestimate}), we will, in this section, assume
that the image of interest is indeed of the form
(\ref{eq1a})--(\ref{eq1b}) and shall denote the associated true
gray values and segmentation by
$\bmu^0=(\mu^0_1,\ldots,\mu_{m^0}^0)$ and
$\bR^0=(R_1^0,\ldots,R_{m^0}^0)$, respectively.

In order to enable large sample results, we impose further technical
conditions. First, to ensure sufficient separation of the regions and
to avoid sets of zero (Lebesgue) measure in the decomposition of
$[0,1]^2$, it
will be assumed throughout that each $R_\nu^0$ contains an open ball
of suitably small radius: for all $\nu=1,\ldots,m^0$, there is
$\bz_\nu\in R_\nu^0$ and $\epsilon>0$ such that
\[
B_\epsilon(\bz_\nu)=\{\bz\in[0,1]^2\colon\|\bz-\bz_\nu\|
<\epsilon\}\subset
R_\nu^0
\]
with \mbox{$\|\cdot\|$} denoting Euclidean norm on $\mathbb{R}^2$. All
candidate segmentations $\bR$ from which the estimate $\hat\bR$ is
produced in any of (\ref{eq3}) to (\ref{eqbicestimate}) are
restricted to satisfy the same condition.

Next, we assume that the set of grid points $\Xi_n$ is dense in
$[0,1]^2$ in the sense that, for all $\epsilon>0$, there is an
$n_0\geq1$ such that
%
\begin{equation}\label{design}
[0,1]^2\subset\bigcup_{i=1}^nB_\epsilon(\bx_i) \qquad\mbox{for all
}n\geq n_0.
\end{equation}

Last, we assume further that the number of grid points in any
given region grows with the sample size (at the same linear rate) and
therefore require that $a_\nu=\lfloor n\alpha_\nu\rfloor$ with
$\sum_\nu\alpha_\nu=1$, where $\lfloor\cdot\rfloor$ denotes the
integer part.

\subsection{Consistency of MDL segmentation}\label{sec3}
We first consider the MDL segmentation solution (\ref{eq3}). Suppose
for now that $m=m^0$ is known, and let $\hat\bR=\argmin_{\bR}\frac
2n\bMDL(m^0,\bR)$. In this case, we have the following strong
consistency result.
%
\begin{theorem}\label{th1}
Let $\{y_i\}$ be the sequence of random variables specified in~(\ref{eqyi}), and assume that $m=m^0$ is known. Then
\[
\hat{\bR}\to\bR^0 \qquad\mbox{with probability one as $n\to\infty$.}
\]
\end{theorem}

The almost sure convergence in the theorem is defined as follows.
Denote by $\prec$ the lexicographical order in $\mathbb{R}^2$, that
is, $a=(a_1,a_2)\prec b=(b_1,b_2)$ if and only if either $a_1<b_1$ or
$a_1=b_1$ and $a_2<b_2$. We assume throughout that any segmentation
$\bR=(R_1,\ldots,R_m)$ satisfies $R_1\prec\cdots\prec R_m$, where~$R_\nu\prec R_\kappa$~if
and only if there is $z_\nu\,{\in}\,R_\nu$ such
that $z_\nu\prec z_\kappa$ for all \mbox{$z_\kappa\,{\in}\,R_\kappa$}. For two
sets~$A$ and $B$, let now $A\Delta B$ be their symmetric difference.
Denote by $\lambda^2$ the~Lebes\-gue\vadjust{\goodbreak} measure in $\mathbb{R}^2$
restricted to $[0,1]^2$ and set $\hat\bR\Delta\bR^0=\bigcup_{\nu
=1}^{m^0}\hat R_\nu\Delta R^0_{\nu}$. Then,\vspace*{1pt} we mean by $\hat{\bR
}\to\bR^0$ with probability one that 
$P(\limsup_n\{\lambda^2(\hat\bR\Delta\bR^0)=0\})=1$. In other
words, the Lebesgue measure of the random sets $\hat\bR\Delta\bR^0$
is zero in the limit with probability one.

The proofs of Theorems \ref{th1} and \ref{th2} below can be
found in the \hyperref[app]{Appendix}.

Of course, in practice, the assumption that $m^0$ is known is
unrealistic. Establishing consistency in the general case of unknown
$m^0$ is, however, substantially more difficult. Even in the simpler
univariate change-point frameworks, where independent variables are
grouped into segments of identical distributions, only special cases
such as normal distributions and exponential families have been
thoroughly investigated; see, for example, \citet{Leecb97} and
\citet{Yao88}. The reason for this is that sharp tail estimates for maxima of
certain squared Gaussian processes are needed which do not hold for
distributions with thicker tails. See Lemma \ref{lem6} below for more details.
Nevertheless, if we assume the noise is normally distributed, we are
able to establish the following consistency result.\vspace*{-2pt}
%
\begin{theorem}\label{th2}
Let $\{y_i\}$ be the sequence of random variables specified
in~(\ref{eqyi}) and assume that the $\{\eps_i\}$ are normally
distributed. Then
\[
\hat m\stackrel{P}{\to}m^0 \qquad\mbox{as } n\to\infty
\]
and
\[
\hat{\bR}\stackrel{P}{\to}\bR^0 \qquad\mbox{as } n\to\infty,
\]
even if the true value $m=m^0$ is unknown. Here $\stackrel{P}{\to}$
indicates convergence in probability.\vspace*{-2pt}
\end{theorem}

The second convergence in probability is defined as follows. Let now
$\hat\bR\Delta\bR^0=\bigcup_{\nu=1}^\mathfrak{m}R_\nu^0\Delta\hat
\bR_\nu$, where $\mathfrak{m}=\min\{m,m^0\}$. Then, in analogy to
the almost sure convergence above, we use the terminology $\hat{\bR
}\stackrel{P}{\to}\bR^0$ to mean that $\lim_nP(\{\lambda^2(\bR
^0\Delta\hat\bR)=0\})=1$. In words, Theorem \ref{th2} asserts
that, if the noise $\eps_i$ is
normal, the MDL method is capable of recovering the true number of
regions as well as the region boundaries as the number of pixels in
the image goes to infinity.\vspace*{-2pt}

\subsection{Consistency of BIC segmentation}
$\!\!\!$The results stated in Theorems~\ref{th1} and \ref{th2} also hold for
the BIC solution given by (\ref{eqbicestimate}). This statement can
be proofed by modifying the proofs for Theorems~\ref{th1}
and \ref{th2}. Details can be found in the \hyperref[app]{Appendix}.\vspace*{-2pt}

\subsection{AIC segmentation is inconsistent}
While being consistent in the special case of known $m=m^0$, the AIC
solution given by (\ref{eqaicestimate}) is, however, inconsistent in
the general case. The main reason is that its penalty term,~$m$, is
independent of the sample size $n$ and does not properly adjust for
the model complexity. Some details are provided in the
\hyperref[app]{Appendix}.\vadjust{\goodbreak}
%
\section{Simulation results}\label{sec4}
Two sets of simulation experiments were conducted to empirically verify
the theoretical results presented above.

\subsection{Experiment 1}
\begin{figure}[t]

\includegraphics{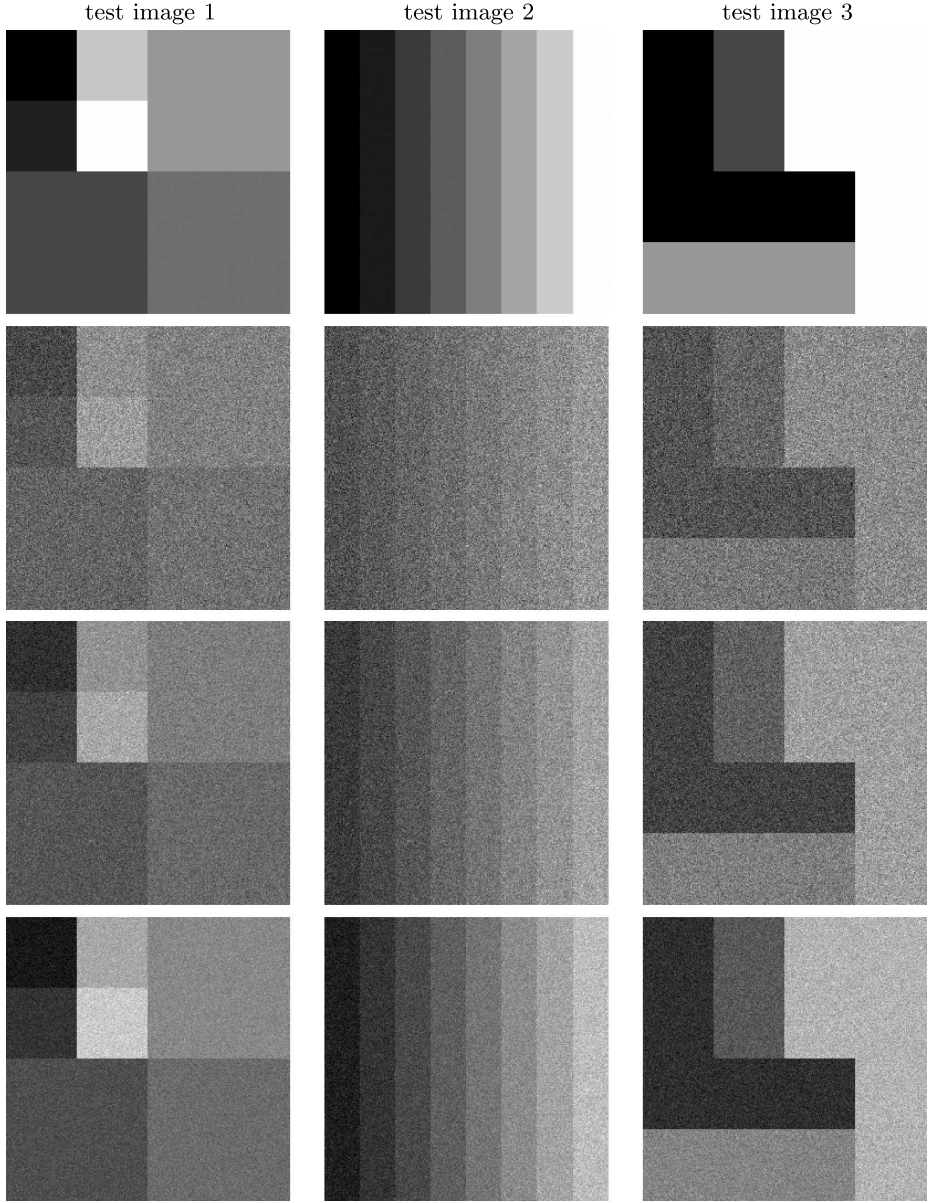}%
\vspace*{-3pt}
\caption{The true test images used in the first numerical experiment
(first row), and typical noisy images generated from $\mbox{snr}=1$ (second
row), 2 (third row) and 4 (last row). All images are plotted with
size $256\times256$.}
\label{figimage}
\vspace*{-3pt}
\end{figure}

Three test images $f$ were used in the first simulation experiment, and they
are displayed in the top row of Figure \ref{figimage}. Recall that
the area and perimeter of each region appear explicitly in the MDL
penalty~(\ref{eq2}), but not the AIC nor the BIC penalty. To assess
the effects of having or not having such quantities as penalty, the
three test images were constructed to have different region areas,
perimeters and area-to-perimeter ratios. Test image 1 has seven
square regions of two different sizes, with true gray values for some
of the adjacent regions being very close. Test image 2 contains eight
rectangular regions of same size, with true gray values increasing from
the left to the right. Test image 3 contains four regions of
different sizes and shapes.

Noisy images were generated by adding Gaussian white noise with
variance $\sigma^2$ to each of the test images. Three signal-to-noise
ratios (snrs) were used: 1, 2 and 4, where snr is defined as
$\sqrt{\var(f)}/\sigma$. Some typical noisy images are also displayed
in Figure \ref{figimage}. Note that for $\mbox{snr}=1$ some of the
region boundaries are hardly visible. Four image sizes were used:
$n=64^2, 128^2$, $256^2$ and $512^2$, and the number of repetitions
for each configuration was~500.

For each noisy image, the AIC, BIC and MDL segmentation
solutions~(\ref{eq3}) to (\ref{eqbicestimate}) were obtained using
the merging algorithm in \citet{Lee00segind}. To verify the result
that $\hat{m} \stackrel{P}{\to} m^0$ (Theorem \ref{th2}), the number
of regions in each segmentation solution was counted and the
corresponding frequencies are tabulated in Tables \ref{tableregion1}
to \ref{tableregion3}. From these tables the following empirical
conclusions can be made:
\begin{itemize}
\item AIC had a strong tendency to over-estimate $m^0$.
\item The performance of BIC improved as $n$ increased, and
occasionally it over-estimated $m^0$.
\item For reasonably large snr and $n$, MDL always correctly estimated
$m^0$.
\item For small snr and $n$, MDL under-estimated $m^0$. As mentioned
before, for such cases some of the region boundaries are hardly
visible (see Figure \ref{figimage}).
\item When comparing the BIC and MDL results, especially from
Table \ref{tableregion3}, it seems that having the region area and
perimeter in the penalty improved the performance.
\end{itemize}

%
\begin{table}
\tabcolsep=0pt
\caption{Frequencies of $\hat{m}$ estimated from the noisy images
generated from test image 1 for different combinations of snr and
$n$. The value of the true $m^0$ is 7}
\label{tableregion1}
\begin{tabular*}{\tablewidth}{@{\extracolsep{\fill}}lrd{3.0}d{3.0}
d{3.0}d{3.0}d{3.0}d{3.0}d{3.0}d{3.0}d{3.0}d{3.0}d{3.0}d{3.0}@{\hspace*{-3pt}}}
\hline
& & \multicolumn{3}{c}{$\bolds{n=64^2}$}
& \multicolumn{3}{c}{$\bolds{n=128^2}$} &
\multicolumn{3}{c}{$\bolds{n=256^2}$} &
\multicolumn{3}{c@{}}{$\bolds{n=512^2}$}
\\[-4pt]
& & \multicolumn{3}{c}{\hrulefill} & \multicolumn{3}{c}{\hrulefill} &
\multicolumn{3}{c}{\hrulefill} & \multicolumn{3}{c@{}}{\hrulefill}
\\
\textbf{snr} & \multicolumn{1}{c}{$\bolds{\hat{m}}$} &
\multicolumn{1}{c}{\textbf{AIC}}
& \multicolumn{1}{c}{\textbf{BIC}} & \multicolumn{1}{c}{\textbf{MDL}}
& \multicolumn{1}{c}{\textbf{AIC}} & \multicolumn{1}{c}{\textbf{BIC}}
& \multicolumn{1}{c}{\textbf{MDL}} & \multicolumn{1}{c}{\textbf{AIC}}
& \multicolumn{1}{c}{\textbf{BIC}} & \multicolumn{1}{c}{\textbf{MDL}}
& \multicolumn{1}{c}{\textbf{AIC}} & \multicolumn{1}{c}{\textbf{BIC}}
& \multicolumn{1}{c@{}}{\textbf{MDL}} \\
\hline
1 & 3\hphantom{$+$} & 0 & 0 & 0 & 0 & 0 & 0 & 0 & 0 & 0 & 0 & 0 & 0 \\
& 4\hphantom{$+$} & 0 & 0 & 0 & 0 & 0 & 0 & 0 & 0 & 0 & 0 & 0 & 0 \\
& 5\hphantom{$+$} & 0 & 0 & 101 & 0 & 0 & 0 & 0 & 0 & 0 & 0 & 0 & 0 \\
& 6\hphantom{$+$} & 0 & 0 & 237 & 0 & 0 & 0 & 0 & 0 & 0 & 0 & 0 & 0 \\
& $\mathbf{7}$\hphantom{$+$} & \mathbf{0} & \mathbf{485} & \mathbf{162} & \mathbf{3} & \mathbf{495} &
\mathbf{500} & \mathbf{6} & \mathbf{499} & \mathbf{500} & \mathbf{0} & \mathbf{500} &
\mathbf{500} \\
& 8\hphantom{$+$} & 18 & 15 & 0 & 10 & 5 & 0 & 15 & 1 & 0 & 14 & 0 & 0 \\
& 9\hphantom{$+$} & 59 & 0 & 0 & 59 & 0 & 0 & 52 & 0 & 0 & 45 & 0 & 0 \\
& 10$+$ & 423 & 0 & 0 & 428 & 0 & 0 & 427 & 0 & 0 & 441 & 0 & 0 \\
[4pt]
2 & 3\hphantom{$+$} & 0 & 0 & 0 & 0 & 0 & 0 & 0 & 0 & 0 & 0 & 0 & 0 \\
& 4\hphantom{$+$} & 0 & 0 & 0 & 0 & 0 & 0 & 0 & 0 & 0 & 0 & 0 & 0 \\
& 5\hphantom{$+$} & 0 & 0 & 0 & 0 & 0 & 0 & 0 & 0 & 0 & 0 & 0 & 0 \\
& 6\hphantom{$+$} & 0 & 0 & 0 & 0 & 0 & 0 & 0 & 0 & 0 & 0 & 0 & 0 \\
& $\mathbf{7}$\hphantom{$+$} & \mathbf{2} & \mathbf{489} & \mathbf{500} & \mathbf{2} & \mathbf{496} &
\mathbf{500} & \mathbf{2} & \mathbf{499} & \mathbf{500} & \mathbf{1} & \mathbf{500} &
\mathbf{500} \\
& 8\hphantom{$+$} & 22 & 11 & 0 & 25 & 4 & 0 & 24 & 1 & 0 & 16 & 0 & 0 \\
& 9\hphantom{$+$} & 63 & 0 & 0 & 79 & 0 & 0 & 65 & 0 & 0 & 52 & 0 & 0 \\
& 10$+$ & 413 & 0 & 0 & 394 & 0 & 0 & 409 & 0 & 0 & 431 & 0 & 0 \\
[4pt]
4 & 3\hphantom{$+$} & 0 & 0 & 0 & 0 & 0 & 0 & 0 & 0 & 0 & 0 & 0 & 0 \\
& 4\hphantom{$+$} & 0 & 0 & 0 & 0 & 0 & 0 & 0 & 0 & 0 & 0 & 0 & 0 \\
& 5\hphantom{$+$} & 0 & 0 & 0 & 0 & 0 & 0 & 0 & 0 & 0 & 0 & 0 & 0 \\
& 6\hphantom{$+$} & 0 & 0 & 0 & 0 & 0 & 0 & 0 & 0 & 0 & 0 & 0 & 0 \\
& $\mathbf{7}$\hphantom{$+$} & \mathbf{3} & \mathbf{487} & \mathbf{500} & \mathbf{0} & \mathbf{498} &
\mathbf{500} & \mathbf{3} & \mathbf{499} & \mathbf{500} & \mathbf{0} & \mathbf{500} &
\mathbf{500} \\
& 8\hphantom{$+$} & 19 & 12 & 0 & 17 & 2 & 0 & 9 & 1 & 0 & 1 & 0 & 0 \\
& 9\hphantom{$+$} & 64 & 0 & 0 & 54 & 0 & 0 & 31 & 0 & 0 & 10 & 0 & 0 \\
& 10$+$ & 414 & 1 & 0 & 429 & 0 & 0 & 457 & 0 & 0 & 489 & 0 & 0 \\
\hline
\end{tabular*}
\end{table}
%

%
\begin{table}
\tabcolsep=0pt
\caption{Similar to Table \protect\ref{tableregion1} but for test
image 2. The value of the true $m^0$ is 8}
\label{tableregion2}
\begin{tabular*}{\tablewidth}{@{\extracolsep{\fill}}lrd{3.0}d{3.0}
d{3.0}d{3.0}d{3.0}d{3.0}d{3.0}d{3.0}d{3.0}d{3.0}d{3.0}d{3.0}@{\hspace*{-3pt}}}
\hline
& & \multicolumn{3}{c}{$\bolds{n=64^2}$}
& \multicolumn{3}{c}{$\bolds{n=128^2}$} &
\multicolumn{3}{c}{$\bolds{n=256^2}$} &
\multicolumn{3}{c@{}}{$\bolds{n=512^2}$}
\\[-4pt]
& & \multicolumn{3}{c}{\hrulefill} & \multicolumn{3}{c}{\hrulefill} &
\multicolumn{3}{c}{\hrulefill} & \multicolumn{3}{c@{}}{\hrulefill}
\\
\textbf{snr} & \multicolumn{1}{c}{$\bolds{\hat{m}}$} &
\multicolumn{1}{c}{\textbf{AIC}}
& \multicolumn{1}{c}{\textbf{BIC}} & \multicolumn{1}{c}{\textbf{MDL}}
& \multicolumn{1}{c}{\textbf{AIC}} & \multicolumn{1}{c}{\textbf{BIC}}
& \multicolumn{1}{c}{\textbf{MDL}} & \multicolumn{1}{c}{\textbf{AIC}}
& \multicolumn{1}{c}{\textbf{BIC}} & \multicolumn{1}{c}{\textbf{MDL}}
& \multicolumn{1}{c}{\textbf{AIC}} & \multicolumn{1}{c}{\textbf{BIC}}
& \multicolumn{1}{c@{}}{\textbf{MDL}} \\
\hline
1 & 3\hphantom{$+$} & 0 & 0 & 213 & 0 & 0 & 0 & 0 & 0 & 0 & 0 & 0 & 0 \\
& 4\hphantom{$+$} & 0 & 0 & 276 & 0 & 0 & 124 & 0 & 0 & 0 & 0 & 0 & 0 \\
& 5\hphantom{$+$} & 0 & 1 & 11 & 0 & 0 & 312 & 0 & 0 & 0 & 0 & 0 & 0 \\
& 6\hphantom{$+$} & 0 & 23 & 0 & 0 & 0 & 57 & 0 & 0 & 0 & 0 & 0 & 0 \\
& 7\hphantom{$+$} & 0 & 127 & 0 & 0 & 0 & 7 & 0 & 0 & 2 & 0 & 0 & 0 \\
& $\mathbf{8}$\hphantom{$+$} & \mathbf{5} & \mathbf{203} & \mathbf{0} & \mathbf{78} & \mathbf{492} &
\mathbf{0} & \mathbf{69} & \mathbf{500} & \mathbf{498} & \mathbf{75} & \mathbf{500} &
\mathbf{500} \\
& 9\hphantom{$+$} & 33 & 114 & 0 & 114 & 6 & 0 & 127 & 0 & 0 & 96 & 0 & 0 \\
& 10$+$ & 462 & 32 & 0 & 308 & 2 & 0 & 304 & 0 & 0 & 329 & 0 & 0 \\
[4pt]
2 & 3\hphantom{$+$} & 0 & 0 & 0 & 0 & 0 & 0 & 0 & 0 & 0 & 0 & 0 & 0 \\
& 4\hphantom{$+$} & 0 & 0 & 12 & 0 & 0 & 0 & 0 & 0 & 0 & 0 & 0 & 0 \\
& 5\hphantom{$+$} & 0 & 0 & 77 & 0 & 0 & 0 & 0 & 0 & 0 & 0 & 0 & 0 \\
& 6\hphantom{$+$} & 0 & 0 & 138 & 0 & 0 & 0 & 0 & 0 & 0 & 0 & 0 & 0 \\
& 7\hphantom{$+$} & 0 & 0 & 147 & 0 & 0 & 0 & 0 & 0 & 0 & 0 & 0 & 0 \\
& $\mathbf{8}$\hphantom{$+$} & \mathbf{85} & \mathbf{488} & \mathbf{126} & \mathbf{82} & \mathbf{500} &
\mathbf{500} & \mathbf{92} & \mathbf{500} & \mathbf{500} & \mathbf{66} & \mathbf{500} &
\mathbf{500} \\
& 9\hphantom{$+$} & 119 & 12 & 0 & 114 & 0 & 0 & 95 & 0 & 0 & 94 & 0 & 0 \\
& 10$+$ & 296 & 0 & 0 & 304 & 0 & 0 & 313 & 0 & 0 & 340 & 0 & 0 \\
[4pt]
4 & 3\hphantom{$+$} & 0 & 0 & 0 & 0 & 0 & 0 & 0 & 0 & 0 & 0 & 0 & 0 \\
& 4\hphantom{$+$} & 0 & 0 & 0 & 0 & 0 & 0 & 0 & 0 & 0 & 0 & 0 & 0 \\
& 5\hphantom{$+$} & 0 & 0 & 0 & 0 & 0 & 0 & 0 & 0 & 0 & 0 & 0 & 0 \\
& 6\hphantom{$+$} & 0 & 0 & 0 & 0 & 0 & 0 & 0 & 0 & 0 & 0 & 0 & 0 \\
& 7\hphantom{$+$} & 0 & 0 & 0 & 0 & 0 & 0 & 0 & 0 & 0 & 0 & 0 & 0 \\
& $\mathbf{8}$\hphantom{$+$} & \mathbf{67} & \mathbf{499} & \mathbf{500} & \mathbf{76} & \mathbf{500} &
\mathbf{500} & \mathbf{84} & \mathbf{500} & \mathbf{500} & \mathbf{65} & \mathbf{500} &
\mathbf{500} \\
& 9\hphantom{$+$} & 96 & 1 & 0 & 126 & 0 & 0 & 102 & 0 & 0 & 115 & 0 & 0 \\
& 10$+$ & 337 & 0 & 0 & 298 & 0 & 0 & 314 & 0 & 0 & 320 & 0 & 0 \\ \hline
\end{tabular*}
\end{table}
%

%
\begin{table}
\tabcolsep=0pt
\caption{Similar to Table \protect\ref{tableregion1} but for test
image 3. The value of the true $m^0$ is 4}
\label{tableregion3}
\begin{tabular*}{\tablewidth}{@{\extracolsep{\fill}}lrd{3.0}d{3.0}
d{3.0}d{3.0}d{3.0}d{3.0}d{3.0}d{3.0}d{3.0}d{3.0}d{3.0}d{3.0}@{\hspace*{-3pt}}}
\hline
& & \multicolumn{3}{c}{$\bolds{n=64^2}$}
& \multicolumn{3}{c}{$\bolds{n=128^2}$} &
\multicolumn{3}{c}{$\bolds{n=256^2}$} &
\multicolumn{3}{c@{}}{$\bolds{n=512^2}$}
\\[-4pt]
& & \multicolumn{3}{c}{\hrulefill} & \multicolumn{3}{c}{\hrulefill} &
\multicolumn{3}{c}{\hrulefill} & \multicolumn{3}{c@{}}{\hrulefill}
\\
\textbf{snr} & \multicolumn{1}{c}{$\bolds{\hat{m}}$} &
\multicolumn{1}{c}{\textbf{AIC}}
& \multicolumn{1}{c}{\textbf{BIC}} & \multicolumn{1}{c}{\textbf{MDL}}
& \multicolumn{1}{c}{\textbf{AIC}} & \multicolumn{1}{c}{\textbf{BIC}}
& \multicolumn{1}{c}{\textbf{MDL}} & \multicolumn{1}{c}{\textbf{AIC}}
& \multicolumn{1}{c}{\textbf{BIC}} & \multicolumn{1}{c}{\textbf{MDL}}
& \multicolumn{1}{c}{\textbf{AIC}} & \multicolumn{1}{c}{\textbf{BIC}}
& \multicolumn{1}{c@{}}{\textbf{MDL}} \\
\hline
1 & 3\hphantom{$+$} & 0 & 0 & 2 & 0 & 0 & 0 & 0 & 0 & 0 & 0 & 0 & 0 \\
& $\mathbf{4}$\hphantom{$+$} & \mathbf{9} & \mathbf{493} & \mathbf{498} & \mathbf{4} & \mathbf{498} &
\mathbf{500} & \mathbf{6} & \mathbf{499} & \mathbf{500} & \mathbf{8} & \mathbf{500} &
\mathbf{500} \\
& 5\hphantom{$+$} & 34 & 6 & 0 & 33 & 2 & 0 & 35 & 1 & 0 & 22 & 0 & 0 \\
& 6\hphantom{$+$} & 63 & 1 & 0 & 60 & 0 & 0 & 70 & 0 & 0 & 57 & 0 & 0 \\
& 7\hphantom{$+$} & 94 & 0 & 0 & 97 & 0 & 0 & 86 & 0 & 0 & 98 & 0 & 0 \\
& 8\hphantom{$+$} & 80 & 0 & 0 & 113 & 0 & 0 & 103 & 0 & 0 & 95 & 0 & 0 \\
& 9\hphantom{$+$} & 99 & 0 & 0 & 96 & 0 & 0 & 87 & 0 & 0 & 88 & 0 & 0 \\
& 10$+$ & 121 & 0 & 0 & 97 & 0 & 0 & 113 & 0 & 0 & 132 & 0 & 0 \\
[4pt]
2 & 3\hphantom{$+$} & 0 & 0 & 0 & 0 & 0 & 0 & 0 & 0 & 0 & 0 & 0 & 0 \\
& $\mathbf{4}$\hphantom{$+$} & \mathbf{6} & \mathbf{494} & \mathbf{500} & \mathbf{7} & \mathbf{498} &
\mathbf{500} & \mathbf{5} & \mathbf{499} & \mathbf{500} & \mathbf{3} & \mathbf{500} &
\mathbf{500} \\
& 5\hphantom{$+$} & 29 & 6 & 0 & 22 & 2 & 0 & 28 & 1 & 0 & 24 & 0 & 0 \\
& 6\hphantom{$+$} & 69 & 0 & 0 & 70 & 0 & 0 & 58 & 0 & 0 & 71 & 0 & 0 \\
& 7\hphantom{$+$} & 92 & 0 & 0 & 97 & 0 & 0 & 87 & 0 & 0 & 85 & 0 & 0 \\
& 8\hphantom{$+$} & 102 & 0 & 0 & 92 & 0 & 0 & 124 & 0 & 0 & 85 & 0 & 0 \\
& 9\hphantom{$+$} & 78 & 0 & 0 & 91 & 0 & 0 & 87 & 0 & 0 & 80 & 0 & 0 \\
& 10$+$ & 124 & 0 & 0 & 121 & 0 & 0 & 111 & 0 & 0 & 152 & 0 & 0 \\
[4pt]
4 & 3\hphantom{$+$} & 0 & 0 & 0 & 0 & 0 & 0 & 0 & 0 & 0 & 0 & 0 & 0 \\
& $\mathbf{4}$\hphantom{$+$} & \mathbf{4} & \mathbf{492} & \mathbf{500} & \mathbf{5} & \mathbf{499} &
\mathbf{500} & \mathbf{3} & \mathbf{500} & \mathbf{500} & \mathbf{2} & \mathbf{500} &
\mathbf{500} \\
& 5\hphantom{$+$} & 29 & 8 & 0 & 24 & 1 & 0 & 12 & 0 & 0 & 4 & 0 & 0 \\
& 6\hphantom{$+$} & 56 & 0 & 0 & 49 & 0 & 0 & 46 & 0 & 0 & 15 & 0 & 0 \\
& 7\hphantom{$+$} & 82 & 0 & 0 & 87 & 0 & 0 & 62 & 0 & 0 & 31 & 0 & 0 \\
& 8\hphantom{$+$} & 102 & 0 & 0 & 104 & 0 & 0 & 101 & 0 & 0 & 44 & 0 & 0 \\
& 9\hphantom{$+$} & 104 & 0 & 0 & 94 & 0 & 0 & 92 & 0 & 0 & 76 & 0 & 0 \\
& 10$+$ & 123 & 0 & 0 & 137 & 0 & 0 & 184 & 0 & 0 & 328 & 0 & 0 \\
\hline
\end{tabular*}
\end{table}

The other major theoretical result that we want to verify is that
$\hat{\bR}$ converges to~$\bR^0$ (Theorems \ref{th1}
and \ref{th2}). However, it is not as straightforward as verifying
$\hat{m} \stackrel{P}{\to} m^0$, as there is no universally agreed
distance metric for measuring the distance between two image
partitions $\hat{\bR}$ and $\bR^0$ [although some related work can be
found in \citet{Baddeley92b}]. To circumvent this issue, we use a
somewhat stricter metric, the mean-squared-error (MSE), defined as
$\mbox{MSE}(\hat{f})=\sum_{i=1}^{n}(f_i-\hat{f}_i)^2$. The reason
we see $\mbox{MSE}(\hat{f})$ as a stricter metric is that, given that
$m^0$ is correctly estimated, it is extremely likely that
$\hat{\bR}=\bR^0$ when $\mbox{MSE}(\hat{f})=0$, but not vice
versa.\vadjust{\goodbreak}

%
\begin{table}
\tabcolsep=2pt
\caption{The averaged $\mbox{MSE}(\hat{f})$ values (multiplied by
1,000) for each combination of test image,\vspace*{1pt} snr and $n$ for the first
simulation experiment. Numbers in parentheses are the ratios $\{\mbox
{MSE}(\hat{f})\}^{0.5}/\sigma$. Boldface indicates the smallest value
for each experimental setting}
\label{tablemse}
{\fontsize{7.5pt}{11pt}\selectfont{
\begin{tabular*}{\tablewidth}{@{\extracolsep{\fill}}lcc@{\hspace*{6pt}}k{2.13}k{3.14}k{1.16}k{1.16}@{}}
\hline
\textbf{Image} & \textbf{snr} & & \multicolumn{1}{c}{$\bolds{n=64^2}$}
& \multicolumn{1}{c}{$\bolds{n=128^2}$} & \multicolumn{1}{c}{$\bolds{n=256^2}$}
& \multicolumn{1}{c@{}}{$\bolds{n=512^2}$} \\
\hline
1 & 1 & AIC & 18,.58\ (0.09352) & 4,.510\ (0.04607) & 1,.090\ (0.02265)
& 0,.2756\ (0.01139) \\
& & BIC & \mathbf{6},\bolds{.}\mathbf{193}\ \mbox{\textbf{(0.05399)}} & 0,.9575\ (0.02123)
& 0,.2244\ (0.01028) & \mathbf{0},\bolds{.}\mathbf{05753}\ \mbox{\textbf{(0.005203)}} \\
& & MDL & 31,.14\ (0.1211) & \mathbf{0},\bolds{.}\mathbf{9304}\ \mbox{\textbf{(0.02092)}}
& \mathbf{0},\bolds{.}\mathbf{2230}\ \mbox{\textbf{(0.01024)}} & \mathbf{0},\bolds{.}\mathbf{05753}\ \mbox{\textbf{(0.005203)}} \\
[4pt]
1 & 2 & AIC & 4,.196\ (0.08887) & 1,.050\ (0.04447) & 0,.2689\ (0.02250) &
0,.06735\ (0.01126) \\
& & BIC & 0,.9305\ (0.04185) & 0,.2291\ (0.02076) & 0,.05672\ (0.01033) &
\mathbf{0},\bolds{.}\mathbf{01441}\ \mbox{\textbf{(0.005208)}} \\
& & MDL & \mathbf{0},\bolds{.}\mathbf{8783}\ \mbox{\textbf{(0.04066)}} & \mathbf{0},\bolds{.}\mathbf{2236}\ \mbox{\textbf{(0.02052)}}
& \mathbf{0},\bolds{.}\mathbf{05630}\ \mbox{\textbf{(0.01029)}} & \mathbf{0},\bolds{.}\mathbf{01441}\ \mbox{\textbf{(0.005208)}} \\
[4pt]
1 & 4 & AIC & 1,.076\ (0.09002) & 0,.2736\ (0.04539) & 0,.06671\ (0.02241) &
0,.01682\ (0.01125) \\
& & BIC & 0,.2472\ (0.04314) & 0,.05934\ (0.02114) & 0,.01424\ (0.01035)
& \mathbf{0},\bolds{.}\mathbf{003550}\ \mbox{\textbf{(0.005170)}} \\
& & MDL & \mathbf{0},\bolds{.}\mathbf{2280}\ \mbox{\textbf{(0.04144)}} & \mathbf{0},\bolds{.}\mathbf{05869}\ \mbox{\textbf{(0.02102)}}
& \mathbf{0},\bolds{.}\mathbf{01414}\ \mbox{\textbf{(0.01032)}} & \mathbf{0},\bolds{.}\mathbf{003550}\ \mbox{\textbf{(0.005170)}} \\
[4pt]
2 & 1 & AIC & \mathbf{76},\bolds{.}\mathbf{23}\ \mbox{\textbf{(0.1894)}} & 6,.908\ (0.05701) & 1,.661\ (0.02796) &
0,.4176\ (0.01402) \\
& & BIC & 112,.8\ (0.2304) & \mathbf{3},\bolds{.}\mathbf{038}\ \mbox{\textbf{(0.03781)}}
& \mathbf{0},\bolds{.}\mathbf{6388}\ \mbox{\textbf{(0.01734)}} & \mathbf{0},\bolds{.}\mathbf{1617}\ \mbox{\textbf{(0.008724)}} \\
& & MDL & 472,.8\ (0.4717) & 218,.2\ (0.3204) & 0,.8846\ (0.02040)
& \mathbf{0},\bolds{.}\mathbf{1617}\ \mbox{\textbf{(0.008724)}} \\
[4pt]
2 & 2 & AIC & 6,.726\ (0.1125) & 1,.655\ (0.05581) & 0,.4015\ (0.02749)
& 0,.1047\ (0.01404) \\
& & BIC & \mathbf{3},\bolds{.}\mathbf{212}\ \mbox{\textbf{(0.07775)}} & \mathbf{0},\bolds{.}\mathbf{6411}\ \mbox{\textbf{(0.03474)}}
& \mathbf{0},\bolds{.}\mathbf{1540}\ \mbox{\textbf{(0.01702)}} & \mathbf{0},\bolds{.}\mathbf{04023}\ \mbox{\textbf{(0.008702)}} \\
& & MDL & 90,.07\ (0.4118) & \mathbf{0},\bolds{.}\mathbf{6411}\ \mbox{\textbf{(0.03474)}}
& \mathbf{0},\bolds{.}\mathbf{1540}\ \mbox{\textbf{(0.01702)}} & \mathbf{0},\bolds{.}\mathbf{04023}\ \mbox{\textbf{(0.008702)}} \\
[4pt]
2 & 4 & AIC & 1,.697\ (0.1130) & 0,.4143\ (0.05585) & 0,.1027\ (0.02780) &
0,.02516\ (0.01376) \\
& & BIC & 0,.6276\ (0.06874) & \mathbf{0},\bolds{.}\mathbf{1552}\ \mbox{\textbf{(0.03419)}}
& \mathbf{0},\bolds{.}\mathbf{03902}\ \mbox{\textbf{(0.01714)}} & \mathbf{0},\bolds{.}\mathbf{01000}\ \mbox{\textbf{(0.008678)}} \\
& & MDL & \mathbf{0},\bolds{.}\mathbf{6248}\ \mbox{\textbf{(0.06859)}} & \mathbf{0},\bolds{.}\mathbf{1552}\ \mbox{\textbf{(0.03419)}}
& \mathbf{0},\bolds{.}\mathbf{03902}\ \mbox{\textbf{(0.01714)}} & \mathbf{0},\bolds{.}\mathbf{01000}\ \mbox{\textbf{(0.008678)}} \\
[4pt]
3 & 1 & AIC & 11,.88\ (0.07476) & 2,.870\ (0.03675) & 0,.7030\ (0.01819)
& 0,.1759\ (0.009098) \\
& & BIC & \mathbf{2},\bolds{.}\mathbf{078}\ \mbox{\textbf{(0.03127)}} & 0,.4024\ (0.01376)
& 0,.09679\ (0.006749) & \mathbf{0},\bolds{.}\mathbf{02367}\ \mbox{\textbf{(0.003338)}} \\
& & MDL & 2,.545\ (0.03461) & \mathbf{0},\bolds{.}\mathbf{3927}\ \mbox{\textbf{(0.01359)}}
& \mathbf{0},\bolds{.}\mathbf{09558}\ \mbox{\textbf{(0.006707)}} & \mathbf{0},\bolds{.}\mathbf{02367}\ \mbox{\textbf{(0.003338)}} \\
[4pt]
3 & 2 & AIC & 2,.932\ (0.07429) & 0,.7225\ (0.03688) & 0,.1822\ (0.01852) &
0,.04568\ (0.009272) \\
& & BIC & 0,.4140\ (0.02792) & 0,.1053\ (0.01408) & 0,.02521\ (0.006889)
& \mathbf{0},\bolds{.}\mathbf{006404}\ \mbox{\textbf{(0.003472)}} \\
& & MDL & \mathbf{0},\bolds{.}\mathbf{3915}\ \mbox{\textbf{(0.02715)}} & \mathbf{0},\bolds{.}\mathbf{1028}\ \mbox{\textbf{(0.01391)}}
& \mathbf{0},\bolds{.}\mathbf{02494}\ \mbox{\textbf{(0.006852)}} & \mathbf{0},\bolds{.}\mathbf{006404}\ \mbox{\textbf{(0.003472)}} \\
[4pt]
3 & 4 & AIC & 0,.7430\ (0.07479) & 0,.1839\ (0.03721) & 0,.04468\ (0.01834) &
0,.01101\ (0.009106) \\
& & BIC & 0,.1106\ (0.02885) & 0,.02441\ (0.01356)
& \mathbf{0},\bolds{.}\mathbf{005919}\ \mbox{\textbf{(0.006676)}} & \mathbf{0},\bolds{.}\mathbf{001478}\ \mbox{\textbf{(0.003336)}} \\
& & MDL & \mathbf{0},\bolds{.}\mathbf{1041}\ \mbox{\textbf{(0.02799)}} & \mathbf{0},\bolds{.}\mathbf{02415}\ \mbox{\textbf{(0.01348)}}
& \mathbf{0},\bolds{.}\mathbf{005919}\ \mbox{\textbf{(0.006676)}} & \mathbf{0},\bolds{.}\mathbf{001478}\ \mbox{\textbf{(0.003336)}} \\
\hline
\end{tabular*}}}
\vspace*{-3pt}
\end{table}

The averaged values of $\mbox{MSE}(\hat{f})$ and
$\{\mbox{MSE}(\hat{f})\}^{0.5}/\sigma$ are listed in
Table \ref{tablemse}, where $\sigma^2$ is the true noise variance.
As expected, the larger the image size~$n$, the smaller these values
are. Also, the corresponding figures from BIC and MDL are
substantially smaller than those from AIC for large $n$. For small~$n$ and snr,
MDL produced poor $\mbox{MSE}(\hat{f})$ values. It is
due to the fact that MDL under-estimates $m^0$.

\subsection{Experiment 2}
\begin{figure}

\includegraphics{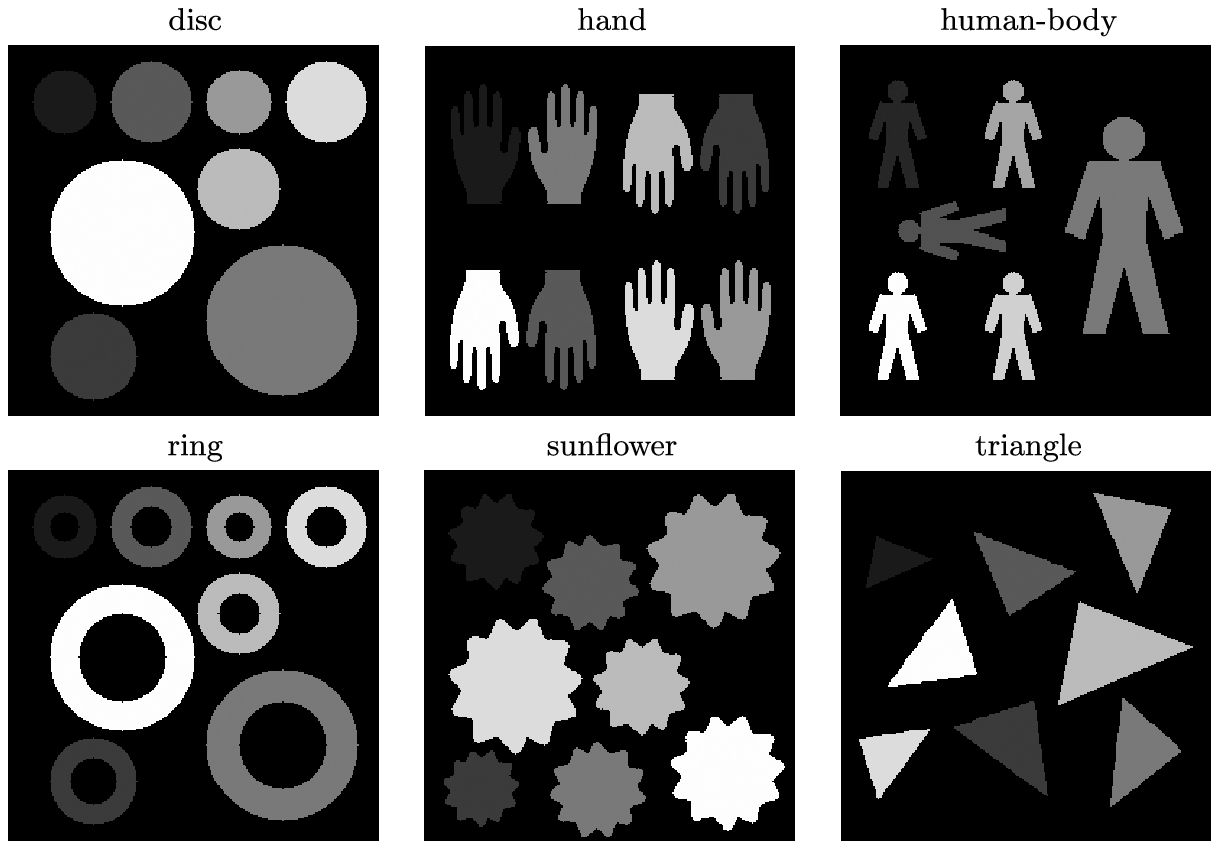}

\caption{The true test images used in the second numerical experiment.}
\label{figimage2}
\end{figure}

Altogether six test images were used in this second numerical
experiment. When comparing to the three test images used in the first
experiments, the shapes of the objects in these six images are more
complicated; see Figure \ref{figimage2}.

We repeated the same testing procedure as above, but only for
$n=256^2$. For each test image, the averages of the estimated number
of regions for AIC, BIC\vadjust{\goodbreak} and MDL segmentation solutions are tabulated
in Table \ref{tablenewmhat}. The standard errors of these averages
are also reported. We have also computed the averaged values of of
$\mbox{MSE}(\hat{f})$ and $\{\mbox{MSE}(\hat{f})\}^{0.5}/\sigma$; they
are listed in Table \ref{tablenewmse}. Empirical conclusions
%
\begin{table}
\caption{The averaged $\hat{m}$ values for the second numerical
experiment. Numbers in parentheses are estimated standard errors. The
true values of $m$ (i.e., $m^0$) are listed in square brackets}
\label{tablenewmhat}
\vspace*{-3pt}
\begin{tabular*}{\tablewidth}{@{\extracolsep{\fill}}lc@{\hspace*{20pt}}k{2.11}k{2.10}k{2.11}@{\hspace*{-0.5pt}}}
\hline
\textbf{Image} $\bolds{[m^0]}$ & & \multicolumn{1}{c}{$\bolds{\mathrm{snr}=1}$}
& \multicolumn{1}{c}{$\bolds{\mathrm{snr}=2}$} & \multicolumn{1}{c@{}}{$\bolds{\mathrm{snr}=4}$} \\ \hline
Disc [8] & AIC & 83,.2\ (0.274) & 69,.0\ (0.268) & 48,.2\ (0.243) \\
& BIC & 20,.9\ (0.165) & 16,.5\ (0.123) & 9,.94\ (0.0689) \\
& MDL & 6,.38\ (0.0219) & 7,.06\ (0.0107) & 8,.05\ (0.014) \\ [4pt]
Hand [8] & AIC & 77,.8\ (0.259) & 63,.7\ (0.247) & 39,.6\ (0.219) \\
& BIC & 20,.4\ (0.139) & 15,.5\ (0.106) & 9,.45\ (0.0636) \\
& MDL & 6,.84\ (0.0259) & 8,.05\ (0.0245) & 8,.13\ (0.0168) \\ [4pt]
Human-body [6] & AIC & 67,.7\ (0.268) & 47,.9\ (0.247) & 25,.3\ (0.187) \\
& BIC & 15,.7\ (0.130) & 8,.97\ (0.0951) & 6,.15\ (0.0194) \\
& MDL & 5,.04\ (0.00964) & 6,.23\ (0.0253) & 6,.03\ (0.00739) \\ [4pt]
Ring [16] & AIC & 81,.1\ (0.266) & 69,.6\ (0.244) & 48,.9\ (0.218) \\
& BIC & 24,.8\ (0.153) & 22,.1\ (0.120) & 16,.7\ (0.0613) \\
& MDL & 11,.2\ (0.0279) & 13,.9\ (0.0184) & 15,.2\ (0.0189) \\ [4pt]
Sunflower [8] & AIC & 81,.8\ (0.289) & 67,.6\ (0.250) & 47,.8\ (0.259) \\
& BIC & 20,.0\ (0.153) & 15,.7\ (0.123) & 10,.2\ (0.0939) \\
& MDL & 6,.07\ (0.0117) & 7,.41\ (0.0222) & 8,.15\ (0.0224) \\ [4pt]
Triangle [8] & AIC & 75,.7\ (0.276) & 62,.6\ (0.248) & 35,.4\ (0.220) \\
& BIC & 18,.6\ (0.138) & 14,.5\ (0.119) & 8,.48\ (0.0313) \\
& MDL & 6,.97\ (0.0101) & 7,.57\ (0.0223) & 7,.99\ (0.00597) \\ \hline
\end{tabular*}
\vspace*{-3pt}
\end{table}
obtainable from these two tables are similar to those from the first
experiment. A~noteworthy observation is that, when snr is not large,
the tendency for BIC to over-estimate $m^0$ is more apparent for these
new test images, that is, when the object boundaries are more
complex.\looseness=-1

\begin{table}
\caption{The\vspace*{0.5pt} averaged $\mbox{MSE}(\hat{f})$ values (multiplied by
1,000) for each combination of test image and snr. Numbers in
parentheses are the ratios
$\{\mbox{MSE}(\hat{f})\}^{0.5}/\sigma$. Boldface indicates
the smallest value for~each experimental setting}
\label{tablenewmse}
\begin{tabular*}{\tablewidth}{@{\extracolsep{\fill}}lc@{\hspace*{27pt}}cck{2.12}@{}}
\hline
\textbf{Image} & & \multicolumn{1}{c}{$\bolds{\mathrm{snr}=1}$}
& \multicolumn{1}{c}{$\bolds{\mathrm{snr}=2}$} & \multicolumn{1}{c@{}}{$\bolds{\mathrm{snr}=4}$} \\ \hline
Disc & AIC & 475.4 (0.2333) & 81.55 (0.1932) & 10,.44\ (0.1383) \\
& BIC & \textbf{405.7 (0.2155)} & \textbf{65.78 (0.1735)}
& 7,.811\ (0.1196) \\
& MDL & 428.7 (0.2215) & 69.56 (0.1784) & \mathbf{7},\mbox{\textbf{.763}}
\ \mbox{\textbf{(0.1192)}} \\
[4pt]
Hand & AIC & 504.9 (0.2950) & 79.51 (0.2342) & 10,.75\ (0.1722) \\
& BIC & \textbf{465.3 (0.2832)} & \textbf{70.62 (0.2207)} & \mathbf{9},\mbox{\textbf{.522}}
\ \mbox{\textbf{(0.1621)}} \\
& MDL & 485.4 (0.2893) & 71.22 (0.2216) & 9,.853\ (0.1649) \\
[4pt]
Human-body & AIC & 135.3 (0.2443) & 19.82 (0.1870) & 1,.491\ (0.1026) \\
& BIC & \textbf{119.9 (0.2300)} & \textbf{17.17 (0.1741)} &
\mathbf{1},\mbox{\textbf{.208}}\ \mbox{\textbf{(0.09234)}} \\
& MDL & 120.9 (0.2309) & 17.53 (0.1759) & 1,.217\ (0.09269) \\
[4pt]
Ring & AIC & 541.1 (0.2774) & 81.89 (0.2158) & 11,.00\ (0.1582) \\
& BIC & \textbf{493.2 (0.2648)} & \textbf{70.89 (0.2008)} &
\mathbf{9},\mbox{\textbf{.314}}\
\mbox{\textbf{(0.1456)}} \\
& MDL & 520.8 (0.2721) & 73.74 (0.2048) & 9,.572\ (0.1476) \\
[4pt]
Sunflower & AIC & 527.3 (0.2517) & 89.43 (0.2073) & 12,.98\ (0.1580) \\
& BIC & \textbf{464.1 (0.2362)} & \textbf{74.97 (0.1898)} &
\mathbf{10},\mbox{\textbf{.54}}\ \mbox{\textbf{(0.1423)}} \\
& MDL & 488.3 (0.2422) & 83.32 (0.2001) & 10,.74\ (0.1437) \\
[4pt]
Triangle & AIC & 219.8 (0.2165) & 32.64 (0.1668) & 3,.242\ (0.1051) \\
& BIC & 182.6 (0.1973) & 24.84 (0.1455) & \mathbf{2},\mbox{\textbf{.326}}
\ \mbox{\textbf{(0.08906)}} \\
& MDL & \textbf{168.9 (0.1897)} & \textbf{23.50 (0.1416)} & 2,.353\ (0.08957)
\\ \hline
\end{tabular*}
\end{table}

\section{Real image segmentation}
Figure \ref{figsegmrSAR}(a) displays a synthetic aperture radar (SAR)
image of a rural area. It is of dimension $250\times250$ and is made
available by Dr. E. Attema of the European Space Research and
Technology Centre. The image has been log-transformed in order to
stabilize the noise variance. It would be useful to segment the image
into regions of similar vegetation.

Notice that the image is extremely noisy (i.e., low snr) and hence
difficult to obtain a good segmentation. Therefore, we applied the
MDL criterion to segment the image, as the simulation results above
suggest\vadjust{\goodbreak} that both AIC and BIC would heavily oversegment the image.
The MDL segmented result, which consists of 34 segmented regions, is
given in Figure \ref{figsegmrSAR}(b).

Even though a Gaussian noise assumption may not be appropriate for
this SAR image, the MDL criterion produced a reasonable segmentation.
The most apparent weakness of the segmentation is the roughness of the
boundaries (many of which should clearly be straight) and the failure
to detect some narrow regions. This weakness can be (at least
partially) attributed to the noisy nature of the image.

\section{Concluding remarks}\label{sec5}
This paper fills an important gap in the image segmentation
literature
by providing a systematic investigation into the theoretical
properties of some popular information theoretic segmentation
methods. It is shown that both the BIC and the MDL segmentation
solutions are\vadjust{\goodbreak} statistically consistent for recovering the number of
objects together with their boundaries in an image. These theoretical
results are empirically verified by simulation experiments. We
also note that our theoretical results can be straightforwardly
extended to higher-dimensional problems, such as volumetric or movie
segmentation.

The numerical results from the simulation experiments also revealed
some discrepancy in the finite sample performances between BIC and
MDL, which can be attributed to the fact that the region area and
perimeter enter explicitly\vadjust{\goodbreak} into the MDL segmentation criterion but not
BIC. These results seem to suggest that, when both the number of
pixels $n$ and the signal-to-noise ratio (snr) are not small, MDL
is capable of producing very stable and reliable results. For those
cases when both $n$ and snr are small, MDL always under-estimated the
number of regions, which led to poor MSE values. However, when one
inspects the noisy images that correspond to such cases, one can see
that, due to the high noise variance, some of the adjacent regions are
hardly distinguishable, which explains the under-estimation of MDL.
Overall the numerical results also suggest that BIC has a tendency to
over-estimate the number of regions, and for those high noise variance
cases, this tendency actually worked in favor of the situation.
Considering all these factors, in practice if the image to be
segmented is not too noisy or not too small in size, one may consider
using MDL, otherwise, use BIC.

\begin{appendix}\label{app}
\section*{Appendix: Proofs}
This Appendix first provides the proofs for Theorems \ref{th1}
and \ref{th2} in Appendices \ref{a1} and \ref{sec42}.
Appendix \ref{a3} covers the BIC and AIC procedures.

\subsection{\texorpdfstring{Proof of Theorem \protect\ref{th1}}{Proof of Theorem 3.1}}\label{a1}
We first provide a number of auxiliary
results and will throughout use the following conventions. The true
segmentation of $[0,1]^2$ will be denoted by
$R_1^0,\ldots,R_{m^0}^0$. All other segmentations will be denoted
$R_1,\ldots,R_m$, while the MDL-based estimates will be $\hat
R_1,\ldots,\hat R_m$. Recall that in the situation of Theorem
\ref{th1}, the number of segments, $m=m^0$, is assumed known.\vadjust{\goodbreak}
%
\begin{lemma}\label{lem1}
Let $y_i=f(\bx_i)+\eps_i$, $i=1,\ldots,n$, be random variables with
$f(\bx)=\mu$ for all $\bx\in[0,1]^2$ and design points $\Xi_n=\{
\bx_1,\ldots,\bx_n\}\subset[0,1]^2$ satisfying (\ref{design}).
Assume furthermore that $\{\eps_i\}$ is a sequence of independent,
identically distributed random variables with zero mean and variance
$\sigma^2$. Fix a subset $R\subset[0,1]^2$, and let $a=\#\mathcal
{A}$ for $\mathcal{A}=\{i\colon\bx_i\in\Xi_n\cap R\}$. Define the
estimators
\[
\hat\mu(R)=\frac1{a}\sum_{i\in\mathcal{A}}y_i  \quad\mbox
{and}\quad
\hat\sigma^2(R)=\frac1a\sum_{i\in\mathcal{A}}\{y_i-\hat\mu(R)\}^2.\vspace*{-2pt}
\]
Then $\hat\mu(R)\to\mu$ and $\hat\sigma^2(R)\to\sigma^2$ with
probability one as $n\to\infty$.\vspace*{-2pt}
\end{lemma}
\begin{pf} Notice that the sequence $\{y_i\}$ is globally independent
and identically distributed with mean $\mu$ and variance $\sigma^2$,
so in particular on any subset $R\subset[0,1]^2$. Both assertions of
the lemma follow therefore directly from the strong law of large
numbers after recognizing that $a\to\infty$ as $n\to\infty$ because
of (\ref{design}).\vspace*{-2pt}
\end{pf}
%
\begin{lemma}\label{lem2}
Let $\{y_i\}$ be the sequence of random variables defined
in~(\ref{eqyi}). Fix a subset $R\subset[0,1]^2$ and denote by
$\hat{\mu}(R)$ the sample mean defined in Lemma \ref{lem1}. Then,
$\hat\mu(R)\to\mu_*$ with probability one, where the limit $\mu_*(R)$
is defined in (\ref{eqnmustar}) below.\vadjust{\goodbreak}
\end{lemma}
\begin{pf} Utilizing the true segmentation, we can write
\[
R=\bigcup_{\nu=1}^mR\cap R_\nu^0=\bigcup_{\ell=1}^2\bigcup_{\nu
\in\mathcal{I}_\ell}R\cap R_\nu^0,
\]
where $\mathcal{I}_1=\{\nu\colon R_\nu^0\subset R\}$ and $\mathcal
{I}_2=\{\nu\colon R\cap R_\nu^0\not=\varnothing\}\setminus\mathcal
{I}_1$, thus ignoring those~$\nu$ for which $R\cap R_\nu^0=\varnothing
$ on the right-hand side of the last display. Define $\tilde a_\nu^0=\#
\tilde{\mathcal{A}}_\nu^0$ for $\tilde{\mathcal{A}}_\nu^0=\{
i\colon\bx_i\in\Xi_n\cap R\cap R_\nu^0\}$ and $a_\nu^0=\#\mathcal
{A}_\nu^0$ for $\mathcal{A}_\nu^0=\{i\colon\bx_i\in\Xi_n\cap
R_\nu^0\}$. It follows from an application of Lemma \ref{lem1} that
%
\begin{eqnarray}\label{eqnmustar}
\hat{\mu}(R)=\frac{1}{a}\sum_{i\in A}y_i
&=&\frac{1}{a}\biggl(\sum_{\nu\in\mathcal{I}_1}\sum_{i\in\mathcal
{A}_\nu^0}y_i
+\sum_{\nu\in\mathcal{I}_2}\sum_{i\in\tilde{\mathcal{A}}_\nu
^0}y_i\biggr) \nonumber\\
&=&\frac{1}{a}\biggl(\sum_{\nu\in\mathcal{ I}_1}a_\nu^0\mu_\nu
^0+\sum_{\nu\in\mathcal{I}_2}\tilde
a_\nu^0\mu_\nu^0\biggr) \\
&\to&\frac{1}{\alpha}\biggl(\sum_{\nu\in\mathcal{I}_1}\alpha_\nu
^0\mu_\nu^0+\sum_{\nu\in\mathcal{I}_2}\tilde\alpha_\nu^0\mu
_\nu^0\biggr)
=:\mu_*(R)
\nonumber
\end{eqnarray}
with probability one as $n\to\infty$, on account of (\ref{design}) and
by assumption on the representation\vspace*{1pt} of the number of design points in
any given region ($a=\lfloor\alpha n\rfloor$, $a_\nu^0=\lfloor
\alpha_\nu^0 n\rfloor$ and $\tilde a_\nu^0=\lfloor
\tilde\alpha_\nu^0n\rfloor$).
\end{pf}

\begin{figure}
\begin{tabular}{@{}cc@{}}

\includegraphics{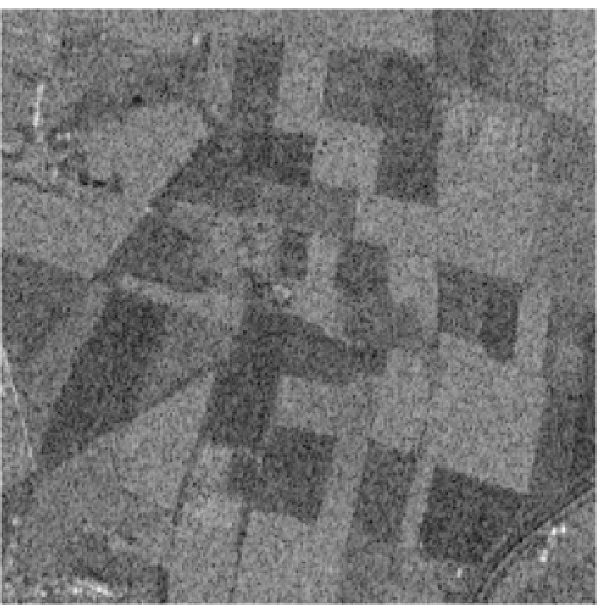}
 & \includegraphics{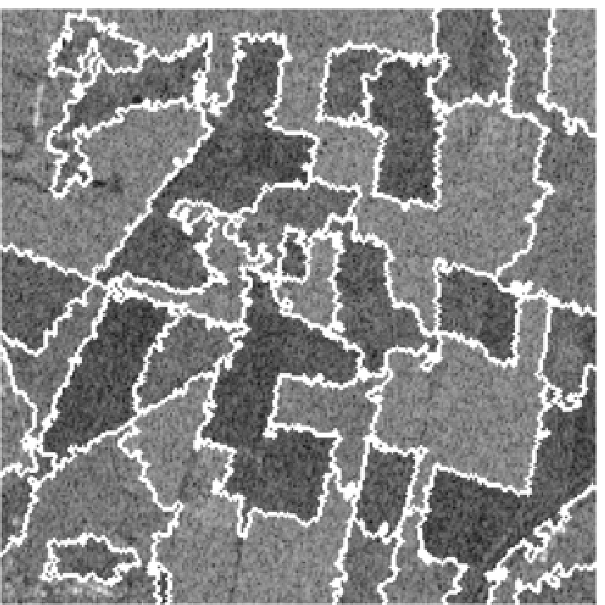}\\
(a) & (b)
\end{tabular}
\caption{Real image segmentation. \textup{(a)}: Observed SAR image
and \textup{(b)}: MDL segmented result.}
\label{figsegmrSAR}
\end{figure}

%
\begin{lemma}\label{lem3}
Let $\{y_i\}$ be the sequence of random variables defined in
(\ref{eqyi}). Fix a subset $R\subset[0,1]^2$ and denote by
$\hat{\sigma}^2(R)$ the variance estimator defined in
Lem\-ma~\ref{lem1}. Then, $\hat\sigma^2(R)\to\sigma^2+\sigma^2_*(R)$
with probability one, where $\sigma^2_*(R)$ is defined in
(\ref{eqnsigmastar}) below.
\end{lemma}
\begin{pf}
Using the notation of the proof of Lemma \ref{lem2} and applying
similar arguments yields the decomposition
\[
\hat{\sigma}^2(R)=\frac1a\sum_{i\in\mathcal{A}}\{y_i-\hat{\mu
}(R)\}^2
=\frac1a\sum_{\nu\in\mathcal{I}_1}\sum_{i\in\mathcal{A}_\nu
^0}\{y_i-\hat{\mu}(R)\}^2
+\frac1a\sum_{\nu\in\mathcal{I}_1}\sum_{i\in\tilde{\mathcal
{A}}_\nu^0}\{y_i-\hat{\mu}(R)\}^2.
\]
Let first $\nu\in\mathcal{I}_1$. By definition of $\mathcal{I}_1$,
$R_\nu^0$ is completely contained in $R$. Therefore, adding and
subtracting the true value $\mu_\nu^0$ from each of the terms
$y_i-\hat\mu(R)$ and subsequently solving the square leads to
\begin{eqnarray*}
\frac1a\sum_{i\in\mathcal{A}_\nu^0} \{y_i-\hat{\mu}(R)\}^2
& = & \frac1a\sum_{i\in\mathcal{A}_\nu^0}(y_i-\mu_\nu^0)^2
-\frac2a\sum_{i\in\mathcal{ A}_\nu^0}(y_i-\mu_\nu^0)\{\mu_\nu
^0-\hat{\mu}(R)\} \\
& &{} +\frac1a\sum_{i\in\mathcal{A}_\nu^0}\{\mu_\nu^0-\hat{\mu
}(R)\}^2 \\
& = & S_{1}+S_{2}+S_{3}.
\end{eqnarray*}
Lemma \ref{lem1} implies for the first term that
\[
S_1=\frac{a_\nu^0}{a}\frac{1}{a_\nu^0}\sum_{i\in\mathcal{A}_\nu
^0}(y_i-\mu_\nu^0)^2\to\frac{\alpha_\nu^0}{\alpha}\sigma^2
\qquad \mbox{a.s.}\qquad (n\to\infty).
\]
The second term $S_2$ is asymptotically small with probability one. To
see this, observe that, by Lemma \ref{lem2}, $\mu_{\nu}^0-\hat{\mu
}(R)$ converges a.s. to $M_\nu^0=\mu_{\nu}^0-\mu_*(R)$ as $n\to
\infty$. For two sequences $\{\xi_n\}$ and $\{\zeta_n\}$ of real
numbers, write $\xi_n\sim\zeta_n$ if $\lim_n\xi_n\zeta_n^{-1}=1$. Then,
using the strong law of large numbers for the i.i.d. sequence~$\{
\eps_i\}$, we obtain that
\[
S_{2}\sim\frac{2M_\nu^0}a\sum_{i\in\mathcal{A}_\nu^0}(y_i-\mu
_\nu^0)
=\frac{2M_\nu^0}a\sum_{i\in\mathcal{A}_\nu^0}\eps_i\to0
\qquad\mbox{a.s.}\qquad (n\to\infty).
\]
Finally, by Lemma \ref{lem2},
\[
S_{3}=\frac{a_\nu^0}{a}\{\mu_{\nu}^0-\hat{\mu}(R)\}^2\to\frac
{\alpha_\nu^0}{\alpha}\{\mu_{\nu}^0-\mu_*(R)\}^2
\qquad \mbox{a.s.}\qquad (n\to\infty).
\]
Let now $\nu\in\mathcal{I}_2$. Then the region $R_\nu^0$ of the
true segmentation is only partially contained in $R$. This means\vspace*{1pt} that,
while all computations can be performed along the blueprint for the
case $\nu\in\mathcal{I}_1$, $\tilde a_\nu^0$, $\alpha_\nu^0$ and
$\tilde{\mathcal{A}}_\nu^0$ have to be used in place of their
respective counterparts $a_\nu^0$, $\alpha_\nu^0$ and $\mathcal
{A}_\nu^0$. Combining these results, we arrive at the almost sure convergence
%
\begin{eqnarray}\label{eqnsigmastar}
\hat\sigma^2(R) & \to&
\frac{\sigma^2}{\alpha}\biggl(\sum_{\nu\in\mathcal{ I}_1}\alpha
_\nu^0+\sum_{\nu\in\mathcal{ I}_2}\tilde\alpha_\nu^0\biggr)
\nonumber\\
& &{}
+\frac1\alpha\biggl[\sum_{\nu\in\mathcal{I}_1}\alpha_\nu^0\{\mu
_\nu^0-\mu_*(R)\}^2
+\sum_{\nu\in\mathcal{ I}_2}\tilde\alpha_\nu^0\{\mu_\nu^0-\mu
_*(R)\}^2\biggr] \\
& = & \sigma^2+\sigma_*^2(R)
\nonumber
\end{eqnarray}
since $\sum_{\mathcal{I}_1}\alpha_\nu^0+\sum_{\mathcal
{I}_2}\tilde\alpha_\nu^0=\alpha$. This proves the assertion.
\end{pf}
%
\begin{lemma}\label{lem4}
Let $\{y_i\}$ be the sequence of random variables defined in~(\ref
{eqyi}). Let $\epsilon>0$ such that, for appropriately chosen\vadjust{\goodbreak} $\bz
_\nu\in R_\nu$ in a segmentation $\bR=(R_1,\ldots,R_m)$,
%
\begin{equation}\label{criterion}
B_\epsilon(\bz_\nu)\subset R_\nu \qquad\mbox{for all }\nu
=1,\ldots,m=m^0.
\end{equation}
Let $\mathcal{R}_\epsilon=\{\bR\colon\bigcup_\nu R_\nu\mbox{
satisfying (\ref{criterion}) such that }a_\nu=\lfloor n\alpha_\nu
\rfloor, \sum_\nu\alpha_\nu=1\}$. Then
\[
\hat\bR=\mathop{\argmin}_{\bR\in\mathcal{R}_\epsilon}\frac2n{
\mathrm{MDL}}(m^0,\bR)\to\bR^0
\qquad \mbox{a.s.}\qquad (n\to\infty),
\]
where $\bR^0$ denotes the true segmentation of\/ $[0,1]^2$.
\end{lemma}
\begin{pf}
Assume that the MDL estimator is not strongly consistent. Thus $\hat
\bR$ does not converge with probability one to\vadjust{\goodbreak} $\bR^0$ as $n\to
\infty$. By boundedness, there exists a monotonically increasing
subsequence $\{n_j\}$ along which $\hat\bR_{n_j}\to\bR^*$ with
probability one, with the limit $\bR^*$ being a member of~$\mathcal
{R}_\epsilon$, and $\lambda^2(\bR^*\Delta\bR^0)>0$ with
probability one. Note that we must have also that $\hat\alpha_\nu\to
\hat\alpha_\nu^*$ along the same subsequence. Note that, with
probability one, $\frac2n{\mathrm{MDL}}(m^0,\bR)\sim\log(\frac
1n\mathrm{RSS}_{m^0})$, where $\sim$ is defined in the proof of Lemma~\ref{lem3}, and that, for $\bR=\bR^*$,
\[
\frac1n\mathrm{RSS}_{m^0}=\frac1n\sum_{\nu=1}^{m^0}\sum_{i\in
\mathcal{A}_\nu^*}\{y_i-\hat\mu(R_\nu^*)\}^2
\]
adopting notation from before. For any $\nu$, there are now two
options: either~$R_\nu^*$ is contained in a region of the true
segmentation, or $R_\nu^*$ has nontrivial intersections with more than
one region of the true segmentation. In the first case, $R_\nu
^*\subset R_\kappa^0$ for some $\kappa$. Hence, Lemma \ref{lem1}
implies that
\[
\frac1n\sum_{i\in\mathcal{A}_\nu^*}\{y_i-\hat\mu(R_\nu^*)\}
^2\to\alpha_\nu^*\sigma^2 \qquad\mbox{a.s.}\qquad (n\to\infty).
\]
In the second case, $R_\nu^*=\bigcup_\kappa R_\nu^*\cap R_\kappa
^0$, where the disjoint union contains at least two elements. Then,
Lemma \ref{lem3} yields that
\[
\frac1n\sum_{i\in\mathcal{A}_\nu^*}\{y_i-\hat\mu(R_\nu^*)\}
^2\to\alpha_\nu^*\sigma^2+\sigma^2_* \qquad\mbox{a.s.}\qquad
(n\to\infty),
\]
where\vspace*{1pt} $\sigma_*^2=\sum_\nu\alpha_\nu^*\sigma^2_*(R_\nu^*)$ with
$\sigma^*(R_\nu^*)$ as in Lemma \ref{lem3}. Observe that, on
account of $\bR^*\not=\bR^0$ [in the sense that $\lambda^2(\bR
^*\Delta\bR^0)\not=0$ almost surely], we have $\sigma_*^2>0$. On
the other hand, $\sigma_*^2=0$ if the true segmentation $\bR^0$ were
used. Consequently, exploiting the continuity and strict concavity of
the logarithm, we arrive at
\begin{eqnarray*}
\lim_{n\to\infty}\frac2n{\mathrm{MDL}}(m^0,\bR^*) & > &
\sum_{\nu=1}^{m^0}\alpha_\nu^0\log\sigma^2=\log\sigma^2
= \lim_{n\to\infty}\frac2n{\mathrm{MDL}}(m^0,\bR^0)\\
&\geq&\lim
_{n\to\infty}\frac2n{\mathrm{MDL}}(m^0,\bR^*),
\end{eqnarray*}
which is a contradiction. Hence, $\hat\bR$ is strongly consistent for
$\bR^0$.
\end{pf}

\subsection{\texorpdfstring{Proof of Theorem \protect\ref{th2}}{Proof of Theorem 3.2}}\label{sec42}

\begin{lemma}\label{lem5}
Let $\{y_i\}$ be the sequence of random variables defined in~(\ref
{eqyi}). If
\[
(\hat m,\hat\bR)=\mathop{\argmin}_{m\leq M\bR\in\mathcal{R}_\epsilon
}\frac2n{\mathrm{MDL}}(m,\bR),
\]
then $P(\hat m\geq m^0)\to1$ as $n\to\infty$.
\end{lemma}
\begin{pf}
Notice that it follows from the proof of Lemma \ref{lem4} that\break $\frac
1n\mathrm{RSS}_{m^0}\to\sigma^2$ with probability one, provided the
true segmentation $\bR^0$ is used in the computations. If $\hat
m<m^0$, then there is at least one $\hat R_\nu$ containing two or more
true regions $R_\kappa^0$. It follows as in the proofs of Lemmas \ref
{lem3} and \ref{lem4} that $P(\frac1n\mathrm{RSS}_m>\sigma
^2+\epsilon)\to1$ as $n\to\infty$ for a suitably chosen $\epsilon
>0$. This implies the claim.
\end{pf}
%
\begin{lemma}\label{lem6}
Let $\{y_i\}$ be the sequence of random variables defined
in~(\ref{eqyi}). If $m^0<m\leq M$, then, for all $\nu=1,\ldots,m^0$,
\[
P\{\hat\bR\in C_\nu^0(n)\}\to0 \qquad(n\to\infty),
\]
where $C_\nu^0(n)=\{\bR=(R_1,\ldots,R_m)\colon\partial R_\kappa
\notin\partial R_\nu^0+B_{\ell(n)}(0), \kappa=1,\ldots,m\}$.
\end{lemma}
\begin{pf}
Fix $1\leq\nu\leq m^0$, and let $\bR\in C^0_\nu(n)$. Because of the
continuity of~$\partial R_\nu^0$, there is a $\bz_\nu\in\partial
R_\nu^0$ such that $\partial R_{\kappa}\cap
B_{\ell(n)}(\bz_\nu)=\varnothing$ for all $\kappa=1,\ldots,m$. Define
$\tilde\bR$ as the segmentation that includes all regions of the form
\[
R_\kappa\cap R_{\nu^\prime}^0\cap B^c_{\ell(n)}(\bz_\nu),\qquad
\kappa=1,\ldots,m; \nu^\prime=1,\ldots,m^0,
\]
and $B_{\ell(n)}(\bz_\nu)$. Clearly,
$\mathrm{RSS}(\bR)\geq\mathrm{RSS}(\tilde\bR)$, where we use the
notations\break $\mathrm{RSS}(\bR)$ and $\mathrm{RSS}(\tilde\bR)$ for
the residual sums of squares based on the respective segmentations $\bR
$ and $\tilde\bR$. Decomposing according
to the true segmentation~$\bR^0$ leads to comparisons of the following
types. Consider first the case $R_{\nu^\prime}^0\cap
B_{\ell(n)}(\bz_\nu)=\varnothing$. Then, it follows as in Lemma 4 of
\citet{Yao88} that
\[
0\leq\sum_{i\in\mathcal{A}_{\nu^\prime}^0}\eps_i^2
-\sum_{\kappa\in\mathcal{I}_{\nu^\prime}}\sum_{i\in\tilde
{\mathcal{A}}_\kappa}\{y_i-\hat\mu(\tilde R_\kappa)\}^2
=\mathcal{O}_P(\ln n)\qquad (n\to\infty),
\]
where $\mathcal{I}_{\nu^\prime}=\{\kappa\colon\tilde R_\kappa
\subset
R_{\nu^\prime}^0\}$,
$\mathcal{A}_{\nu^\prime}^0=\{i\colon\bx_i\in\Xi_n\cap
R_{\nu^\prime}^0\}$ and
$\tilde{\mathcal{A}}_{\kappa}=\{i\colon\bx_i\in\Xi_n\cap\tilde
R_{\kappa}\}$. The rate on the right-hand side of the last display
explicitly uses that the noise $\{\eps_i\}$ follows a normal law and
does not need to be true for arbitrary noise distributions [compare the
remark on page 188 of \citet{Yao88}]. Consider next the case
$R_{\nu^\prime}^0\cap B_{\ell(n)}(\bz_\nu)\not=\varnothing$.
Observe that
the number of design points in $B_{\ell(n)}(\bz_\nu)$ is proportional
to $\ln^2n$, while the number of design points in any $\tilde R_\nu$ is
proportional to the sample size $n$. Any region $\tilde
R_\nu\in\tilde\bR$ obtained from a nontrivial intersection with
$B_{\ell(n)}^c(\bz_\nu)$ has therefore the number of elements reduced
by a factor proportional to $\ln^2 n$. This, however, is negligible
compared to $n$ in the long run. Therefore, the same arguments as
before imply also that
\[
0\leq\sum_{i\in\mathcal{C}^0_{\nu^\prime}}\eps_i^2
-\sum_{\kappa\in\mathcal{J}_{\nu^\prime}}\sum_{i\in\tilde
{\mathcal{A}}_\kappa}\{y_i-\hat\mu(\tilde R_\kappa)\}^2
=\mathcal{O}_P(\ln n)\qquad (n\to\infty),
\]
where $\mathcal{C}_{\nu^\prime}^0=\mathcal{A}_{\nu^\prime
}^0\setminus\mathcal{ B}_{\nu^\prime}$ with $\mathcal{ B}_{\nu
^\prime}=\{i\colon\bx_i\in\Xi_n\cap B_{\ell(n)}(\bz_\nu)\cap
R_{\nu^\prime}^0\}$, and $\mathcal{J}_{\nu^\prime}=\{\kappa\colon
\tilde
R_\kappa\subset R_{\nu^\prime}^0\cap B_{\ell(n)}^c(\bz_\nu)\}$. It
remains\vadjust{\goodbreak} to investigate the region $B_{\ell(n)}(\bz_\nu)$
itself. Without loss of generality assume that $B_{\ell(n)}(\bz_\nu)$
intersects, apart from $R_\nu^0$, only one more true regions
$R_{\nu^\prime}^0$ as the general case can be handled in a similar
fashion. Notice that $b=\# \{B_{\ell(n)}(\bz_\nu)\cap\Xi_n\}
=\lfloor\beta n\rfloor\sim\ln^2n$ by definition. Let furthermore
$b_\nu=\# \{\Xi_n\cap R_\nu^0\cap B_{\ell(n)}(\bz_\nu)\}$ and
$b_{\nu^\prime}=\#\{\Xi_n\cap R_{\nu^\prime}^0\cap B_{\ell
(n)}(\bz_\nu)\}$. Then, we must have $b_\nu=\lfloor\beta_\nu
n\rfloor\sim\ln^2n$ and $b_{\nu^\prime}=\lfloor\beta_{\nu
^\prime}n\rfloor\sim\ln^2n$ for appropriate $\beta_\nu$ and
$\beta_{\nu^\prime}$ satisfying $\beta_\nu+\beta_{\nu^\prime
}=\beta$. Now, utilizing that $y_i-\hat\mu(B_{\ell(n)}(\bz_\nu
))=\eps_i+\mu_\nu-\hat\mu(B_{\ell(n)}(\bz_\nu))$ on $R_\nu^0$
and $y_i-\hat\mu(B_{\ell(n)}(\bz_\nu))=\eps_i+\mu_{\nu^\prime
}-\hat\mu(B_{\ell(n)}(\bz_\nu))$ on $R_{\nu^\prime}^0$, we
obtain that
\begin{eqnarray*}
&&\frac1b\biggl[\sum_{i\in\mathcal{B}_\nu^*}\epsilon_i^2-\sum
_{i\in\mathcal{B}_\nu^*}\bigl\{y_i-\hat\mu\bigl(B_{\ell(n)}(\bz_\nu)\bigr)\bigr\}
^2\biggr]\\
&&\qquad =\frac1b\bigl[b_\nu\bigl\{\mu_\nu-\hat\mu\bigl(B_{\ell(n)}(\bz
_\nu)\bigr)\bigr\}^2+b_{\nu^\prime}\bigl\{\mu_{\nu^\prime}-\hat\mu\bigl(B_{\ell
(n)}(\bz_\nu)\bigr)\bigr\}^2\bigr]+o(1)\\
&&\qquad \to-\frac{\beta_\nu\beta_{\nu^\prime}}{\beta^2}(\mu_\nu
-\mu_{\nu^\prime})^2=B
\end{eqnarray*}
with probability one as $n\to\infty$, where $\mathcal{B}_\nu^*=\{
i\colon\bx_i\in\Xi_n\cap B_{\ell(n)}(\bz_\nu)\}$ and the limit
is clearly negative. Combining the results in the last three displays,
we arrive consequently at
\[
\frac1b\{\mathrm{RSS}-\mathrm{RSS}(\tilde\bR)\}
\stackrel{R}{\to} B<0,
\]
where $\mathrm{RSS}=\sum_{i=1}^n\eps_i^2$. Thus,
\[
\lim_{n\to\infty}\min_{\bR\in[C^0_\nu(n)]^c}\mathrm{RSS}(\bR)>%
\lim_{n\to\infty}\mathrm{RSS}
\geq\lim_{n\to\infty}\mathrm{RSS}(\hat\bR)
\]
with probability approaching one. This implies the assertion.
\end{pf}
%
\begin{lemma}\label{lem7}
Let $\{y_i\}$ be the sequence of random variables defined in~(\ref{eqyi}). If $m^0<m\leq M$ and $\epsilon>0$, then
\[
P\{\mathrm{RSS}-\mathrm{RSS}(\hat\bR)\in[0,L_n(\epsilon
,\hat\bR)]\}\to1 \qquad (n\to\infty),
\]
where $\mathrm{RSS}=\sum_{i=1}^n\eps_i^2$, $\mathrm{RSS}(\hat\bR
)$ is
the residual sum of squares based on the segmentation $\hat\bR=(\hat
R_1,\ldots,\hat R_m)$ selected by the MDL criterion and
$L_n(\epsilon,\bR)=\sigma^2\{\epsilon+2(m-m^0-1)(1+\epsilon)\}\ln n$.
\end{lemma}
\begin{pf}
It follows from Lemma \ref{lem6} that $\hat\bR\in B^0(n)=\bigcap
_{\nu=1}^{m^0}[C_\nu^0(n)]^c$ with probability approaching one. It is
therefore sufficient to verify the claim for an arbitrary segmentation
$\bR\in B^0(n)$. Given such an $\bR$ introduce the finer $\tilde\bR
$ as the segmentation containing the regions
%
\begin{equation}\label{R1}
R_\kappa\cap R_{\nu^\prime}^0\cap[B^0(n)]^c,\qquad \kappa=1,\ldots
,m; \nu^\prime=1,\ldots,m^0,
\end{equation}
and
%
\begin{equation}\label{R2}
R_\kappa\cap R_{\nu^\prime}^0\cap B_\nu^0(n),\qquad \kappa
=1,\ldots,m; \nu,\nu^\prime=1,\ldots,m^0.
\end{equation}
Denote\vspace*{1pt} the collection of regions (\ref{R1}) by $\tilde\bR_1$ and the
collection of regions (\ref{R2}) by~$\tilde\bR_2$. We then have
$\mathrm{RSS}\geq\mathrm{RSS}(\bR)\geq\mathrm{RSS}(\tilde\bR
)=\mathrm{RSS}(\tilde\bR_1)+\mathrm{RSS}(\tilde\bR_2)$. The
number of design points in $\tilde\bR_2$ is, by definition of the
sets $C_\nu^0(n)$, proportional to $\ln n$. An application of Lemma 1
in \citet{Yao88} yields therefore that
\[
\biggl|\sum_{\tilde R_\nu\in\tilde\bR_2}\sum_{i\in{\tilde A}_\nu
}\eps_i^2-\mathrm{RSS}(\tilde\bR_2)\biggr|=\mathcal{O}_P(\ln\ln n)
\qquad (n\to\infty).
\]
For $\tilde R_\nu\in\tilde\bR_1$, let $\tilde a_\nu=\#\tilde R_\nu
$. Since $\bR\in C^0(n)$, it holds that $\#\tilde\bR_1\leq m-m^0$.
As in (17)--(19) of \citet{Yao88}, we conclude therefore with Theorem~2 of \citet{Darling-Erdos56} that, for any $\epsilon>0$ and with
probability approaching one,
\[
\sum_{\tilde R_\nu\in\tilde\bR_1}\sum_{i\in{\tilde A}_\nu}\eps_i^2
\geq\mathrm{RSS}(\tilde\bR_1)
\geq\sum_{\tilde R_\nu\in\tilde\bR_1}\sum_{i\in{\tilde A}_\nu
}\eps_i^2-L_n(\epsilon,\bR).
\]
This completes the proof.
\end{pf}
%
\begin{lemma}\label{lem8}
Let $\{y_i\}$ be the sequence of random variables defined in~%
(\ref{eqyi}). If $m>m^0$, then using the notation of (\ref{eq2}),
it holds for the penalty terms arising from the area and the perimeter
pieces that
\[
\sum_{\kappa=1}^m\ln a_\kappa-\sum_{\nu=1}^{m^0}\ln a_\nu^0\geq
0 \quad\mbox{and}\quad
\sum_{\kappa=1}^mb_\kappa-\sum_{\nu=1}^{m^0}b_\nu^0\geq0
\]
with probability approaching one as $n\to\infty$.
\end{lemma}
\begin{pf}
Lemma \ref{lem6} implies that the oversegmentation $\hat\bR_m$
approximates the true segmentation $\bR^0$ in the sense that, with
probability approaching one, each perimeter $\partial R^0_\nu$ is
uniformly approximated by one or more perimeters $\partial\hat
R_\kappa$. This yields in particular that, for a suitable $\nu_\kappa
=1,\ldots,m^0$, $P(\hat R_\kappa\subset R_{\nu_\kappa}^0)\to1$ for
all $\kappa=1,\ldots,m$. By assumption, we can write that $a_\kappa
=\lambda_{\kappa,\nu}a_{\nu_\kappa}^0$ with $\lambda_{\kappa,\nu
}\to\alpha_\kappa/\alpha_{\nu_\kappa}^0$ as $n\to\infty$. Let
$\mathcal{V}_\nu=\{\kappa^\prime\colon R_{\kappa^\prime}\cap
R_\nu^0\not=\varnothing\}$. Then, with probability approaching one,
\[
\prod_{\kappa=1}^ma_\kappa\Biggl[\prod_{\nu=1}^{m^0}a_\nu^0\Biggr]^{-1}
=\prod_{\nu=1}^{m^0}\prod_{\kappa\in\mathcal{V}_\nu}\lambda
_{\kappa,\nu}(a_\nu^0)^{\#\mathcal{V}_\nu-1}
\geq(\min a_\nu^0)^{m-m^0}\prod_{\nu=1}^{m^0}\prod_{\kappa\in
\mathcal{V}_\nu}\lambda_{\kappa,\nu}
\geq1
\]
since $\sum_\nu(\#\mathcal{V}_\nu-1)=m-m^0$, $a_\nu^0=\lfloor
\alpha_\nu^0n\rfloor$ and the product over the $\lambda_{\kappa
,\nu}$ converges to a finite limit as $n\to\infty$. This implies the
first statement of the lemma. The second claim follows along similar
lines from the fact that the true segmentation ``shares'' all its
perimeters with the oversegmentation with probability approaching one.
Since $m>m^0$, there must at least be one additional perimeter piece
and the assertion follows.
\end{pf}
%
\begin{lemma}\label{lem9}
Let $\{y_i\}$ be the sequence of random variables defined in~(\ref{eqyi}). If $m>m^0$, then
\[
\Delta(m,m^0)=
\frac n2\biggl\{\ln\biggl(\frac{\mathrm{RSS}_m}n\biggr)-\ln
\biggl(\frac{\mathrm{RSS}_{m^0}}n\biggr)\biggr\}+(m-m^0)\ln n\geq0
\]
with probability approaching one as $n\to\infty$.
\end{lemma}
\begin{pf}
Let $\epsilon>0$. By the law of large numbers, we have that $\mathrm
{RSS}=\sum_{i=1}^n\eps_i^2>n(\sigma^2-\epsilon)$. Also, $\mathrm
{RSS}\geq\mathrm{RSS}_{m^0}$. Hence,
\begin{eqnarray*}
\Delta(m,m^0)
&\geq&\frac n2\biggl\{\ln\biggl(\frac{\mathrm{RSS}_m}n\biggr)-\ln
\biggl(\frac{\mathrm{RSS}}n\biggr)\biggr\}+(m-m^0)\ln n \\
&=&\frac n2\ln\biggl(1-\frac{\mathrm{RSS}-\mathrm{RSS}_m}{\mathrm
{RSS}}\biggr)+(m-m^0)\ln n \\
&\geq&\frac n2\ln\biggl\{1-\frac{L_n(\epsilon,\hat\bR)}{n(\sigma
^2-\epsilon)}\biggr\}+(m-m^0)\ln n,
\end{eqnarray*}
where the last inequality follows after an application of Lemma \ref
{lem7}. Continuing as in \citet{Yao88}, using the fact that $\ln
(1-x)>-x(1+\epsilon)$ for small positive $x$ and the definition of
$L_n(\epsilon,\hat\bR)$, the right-hand side can be estimated from
below by
%
\begin{equation}\label{last}
-\frac{\sigma^2(1+\epsilon)}{2(\sigma^2-\epsilon)}\{\epsilon
+2(m-m^0-1)(1+\epsilon)\}\ln
n+(m-m^0)\ln n,
\end{equation}
which is positive with probability approaching one whenever $\epsilon$
is sufficiently small.
\end{pf}

This implies that $\hat m\stackrel{P}{\to} m^0$. The second claim of
Theorem \ref{th2} follows from $P(\mathcal{L}_n)\geq P(\mathcal
{L}_n,\hat m=m^0)\to1$, where $\mathcal{L}_n=\{\lambda^2(\bR
^0\Delta\hat\bR)=0\}$.

\subsection{Proofs for BIC and AIC segmentations}\label{a3}

The counterparts of Theorem~\ref{th1} for the AIC and BIC procedures
are verbatim the same as for the MDL procedure. Consistency in the
case of known $m=m^0$ does therefore not depend on the particular
penalty terms.

The situation is, however, very different in the general case of an
unknown number of segments in the partition. Here, we can prove the
consistency result of Theorem \ref{th2} only for the BIC
procedure. Following the lines of the proofs in Appendix~\ref{sec42},
it can be seen that Lemmas \ref{lem5}--\ref{lem7} deal only with the
RSS term and hold irrespective of the specific penalty term. Lemma~\ref{lem8} deals with the complexity of areas and perimeters unique
to the MDL criterion. The crucial point is therefore Lemma
\ref{lem9}. Repeating\vspace*{1pt} the arguments in its proof, one can for the BIC
criterion similarly verify that, if $m>m^0$,
\[
\tilde\Delta(m,m^0)=
\frac n2\biggl\{\ln\biggl(\frac{\mathrm{RSS}_m}n\biggr)-\ln
\biggl(\frac{\mathrm{RSS}_{m^0}}n\biggr)\biggr\}+(m-m^0)\ln n\geq0
\]
with probability approaching one as $n\to\infty$, utilizing
\[
-\frac{\sigma^2(1+\epsilon)}{2(\sigma^2-\epsilon)}\{\epsilon
+2(m-m^0-1)(1+\epsilon)\}\ln
n+(m-m^0)\ln n
\]
instead of (\ref{last}). This implies consistency of the BIC
procedure. For the AIC segmentation, however, the second term in the
last display becomes $2(m-m^0)$ which grows too slowly to ensure
positivity. Hence AIC-based procedures are inconsistent if $m$ is
unknown.
\end{appendix}

\section*{Acknowledgments}

The authors are grateful to the reviewers and
the Associate Editor for their most useful comments.



\printaddresses


\begin{thebibliography}{22}

\bibitem[\protect\citeauthoryear{Akaike}{1974}]{Akaike74}
\begin{barticle}[mr]
\bauthor{\bsnm{Akaike},~\bfnm{Hirotugu}\binits{H.}}
(\byear{1974}).
\btitle{A new look at the statistical model identification}.
\bjournal{IEEE Trans. Automat. Control}
\bvolume{AC-19}
\bpages{716--723}.
\bnote{System identification and time-series analysis}.
\bid{issn={0018-9286}, mr={0423716}}
\bptok{imsref}%
\end{barticle}
\endbibitem

\bibitem[\protect\citeauthoryear{Baddeley}{1992}]{Baddeley92b}
\begin{barticle}[mr]
\bauthor{\bsnm{Baddeley},~\bfnm{A.~J.}\binits{A.~J.}}
(\byear{1992}).
\btitle{Errors in binary images and an {$L\sp p$} version of the {H}ausdorff
  metric}.
\bjournal{Nieuw Arch. Wisk. (4)}
\bvolume{10}
\bpages{157--183}.
\bid{issn={0028-9825}, mr={1218662}}
\bptok{imsref}%
\end{barticle}
\endbibitem

\bibitem[\protect\citeauthoryear{Darling and Erd{\"o}s}{1956}]{Darling-Erdos56}
\begin{barticle}[mr]
\bauthor{\bsnm{Darling},~\bfnm{D.~A.}\binits{D.~A.}} \AND
  \bauthor{\bsnm{Erd{\"o}s},~\bfnm{P.}\binits{P.}}
(\byear{1956}).
\btitle{A limit theorem for the maximum of normalized sums of independent
  random variables}.
\bjournal{Duke Math. J.}
\bvolume{23}
\bpages{143--155}.
\bid{issn={0012-7094}, mr={0074712}}
\bptok{imsref}%
\end{barticle}
\endbibitem

\bibitem[\protect\citeauthoryear{Glasbey and Horgan}{1995}]{Glasbey-Horgan95}
\begin{bbook}[author]
\bauthor{\bsnm{Glasbey},~\bfnm{Chris~A.}\binits{C.~A.}} \AND
  \bauthor{\bsnm{Horgan},~\bfnm{Graham~W.}\binits{G.~W.}}
(\byear{1995}).
\btitle{Image Analysis for the Biological Sciences}.
\bpublisher{Wiley}, \baddress{Chichester, New York}.
\bptok{imsref}%
\end{bbook}
\endbibitem

\bibitem[\protect\citeauthoryear{Haralick and
  Shapiro}{1992}]{Haralick-Shapiro92}
\begin{bbook}[author]
\bauthor{\bsnm{Haralick},~\bfnm{Robert~M.}\binits{R.~M.}} \AND
  \bauthor{\bsnm{Shapiro},~\bfnm{Linda~G.}\binits{L.~G.}}
(\byear{1992}).
\btitle{Computer and Robot Vision}.
\bpublisher{Addison-Wesley}, \baddress{Reading, MA}.
\bptok{imsref}%
\end{bbook}
\endbibitem

\bibitem[\protect\citeauthoryear{Kanungo et~al.}{1995}]{Kanungo-et-al95}
\begin{bmisc}[author]
\bauthor{\bsnm{Kanungo},~\bfnm{Tapas}\binits{T.}},
  \bauthor{\bsnm{Dom},~\bfnm{Byron}\binits{B.}},
  \bauthor{\bsnm{Niblack},~\bfnm{Wayne}\binits{W.}},
  \bauthor{\bsnm{Steele},~\bfnm{David}\binits{D.}} \AND
  \bauthor{\bsnm{Sheinvald},~\bfnm{Jacob}\binits{J.}}
(\byear{1995}).
\bhowpublished{MDL-based multi-band image segmentation using a fast region merging
scheme.
Technical Report RJ 9960 (87919),
IBM Research Division.}
\bptok{imsref}%
\end{bmisc}
\endbibitem

\bibitem[\protect\citeauthoryear{LaValle and
  Hutchinson}{1995}]{LaValle-Hutchinson95}
\begin{barticle}[author]
\bauthor{\bsnm{LaValle},~\bfnm{Steven~M.}\binits{S.~M.}} \AND
  \bauthor{\bsnm{Hutchinson},~\bfnm{Seth~A.}\binits{S.~A.}}
(\byear{1995}).
\btitle{A {B}ayesian segmentation methodology for parametric image models}.
\bjournal{IEEE Transactions on Pattern Analysis and Machine Intelligence}
\bvolume{17}
\bpages{211--217}.
\bptok{imsref}%
\end{barticle}
\endbibitem

\bibitem[\protect\citeauthoryear{Leclerc}{1989}]{Leclerc89}
\begin{barticle}[author]
\bauthor{\bsnm{Leclerc},~\bfnm{Yvan~G.}\binits{Y.~G.}}
(\byear{1989}).
\btitle{Constructing simple stable descriptions for image partitioning}.
\bjournal{Int. J. Comput. Vis.}
\bvolume{3}
\bpages{73--102}.
\bptok{imsref}%
\end{barticle}
\endbibitem

\bibitem[\protect\citeauthoryear{Lee}{1997}]{Leecb97}
\begin{barticle}[mr]
\bauthor{\bsnm{Lee},~\bfnm{Chung-Bow}\binits{C.-B.}}
(\byear{1997}).
\btitle{Estimating the number of change points in exponential families
  distributions}.
\bjournal{Scand. J. Stat.}
\bvolume{24}
\bpages{201--210}.
\bid{doi={10.1111/1467-9469.t01-1-00058}, issn={0303-6898}, mr={1455867}}
\bptok{imsref}%
\end{barticle}
\endbibitem

\bibitem[\protect\citeauthoryear{Lee}{1998}]{Lee98segcor}
\begin{barticle}[author]
\bauthor{\bsnm{Lee},~\bfnm{Thomas C.~M.}\binits{T.~C.~M.}}
(\byear{1998}).
\btitle{Segmenting images corrupted by correlated noise}.
\bjournal{IEEE Transactions on Pattern Analysis and Machine Intelligence}
\bvolume{20}
\bpages{481--492}.
\bptok{imsref}%
\end{barticle}
\endbibitem

\bibitem[\protect\citeauthoryear{Lee}{2000}]{Lee00segind}
\begin{barticle}[mr]
\bauthor{\bsnm{Lee},~\bfnm{Thomas C.~M.}\binits{T.~C.~M.}}
(\byear{2000}).
\btitle{A minimum description length-based image segmentation procedure, and
  its comparison with a cross-validation-based segmentation procedure}.
\bjournal{J.~Amer. Statist. Assoc.}
\bvolume{95}
\bpages{259--270}.
\bid{issn={0162-1459}, mr={1803154}}
\bptok{imsref}%
\end{barticle}
\endbibitem

\bibitem[\protect\citeauthoryear{Luo and
  Khoshgoftaar}{2006}]{Luo-Khoshgoftaar06}
\begin{barticle}[pbm]
\bauthor{\bsnm{Luo},~\bfnm{Qiming}\binits{Q.}} \AND
  \bauthor{\bsnm{Khoshgoftaar},~\bfnm{Taghi~M.}\binits{T.~M.}}
(\byear{2006}).
\btitle{Unsupervised multiscale color image segmentation based on MDL
  principle}.
\bjournal{IEEE Trans. Image Process.}
\bvolume{15}
\bpages{2755--2761}.
\bid{issn={1057-7149}, pmid={16948319}}
\bptok{imsref}%
\end{barticle}
\endbibitem

\bibitem[\protect\citeauthoryear{Murtagh, Raftery and
  Starck}{2005}]{Murtagh-et-al05}
\begin{barticle}[author]
\bauthor{\bsnm{Murtagh},~\bfnm{F.}\binits{F.}},
  \bauthor{\bsnm{Raftery},~\bfnm{A.~E.}\binits{A.~E.}} \AND
  \bauthor{\bsnm{Starck},~\bfnm{J.~L.}\binits{J.~L.}}
(\byear{2005}).
\btitle{Bayesian inference for multiband image segmentation via model-based
  cluster trees}.
\bjournal{Image and Vision Computing}
\bvolume{23}
\bpages{587--596}.
\bptok{imsref}%
\end{barticle}
\endbibitem

\bibitem[\protect\citeauthoryear{Rissanen}{1989}]{Rissanen89}
\begin{bbook}[mr]
\bauthor{\bsnm{Rissanen},~\bfnm{Jorma}\binits{J.}}
(\byear{1989}).
\btitle{Stochastic Complexity in Statistical Inquiry}.
\bseries{World Scientific Series in Computer Science}
\bvolume{15}.
\bpublisher{World Scientific}, \baddress{Teaneck, NJ}.
\bid{mr={1082556}}
\bptok{imsref}%
\end{bbook}
\endbibitem

\bibitem[\protect\citeauthoryear{Rissanen}{2007}]{Rissanen07}
\begin{bbook}[mr]
\bauthor{\bsnm{Rissanen},~\bfnm{Jorma}\binits{J.}}
(\byear{2007}).
\btitle{Information and Complexity in Statistical Modeling}.
\bpublisher{Springer}, \baddress{New York}.
\bid{mr={2287233}}
\bptok{imsref}%
\end{bbook}
\endbibitem

\bibitem[\protect\citeauthoryear{Schwarz}{1978}]{Schwarz78}
\begin{barticle}[mr]
\bauthor{\bsnm{Schwarz},~\bfnm{Gideon}\binits{G.}}
(\byear{1978}).
\btitle{Estimating the dimension of a model}.
\bjournal{Ann. Statist.}
\bvolume{6}
\bpages{461--464}.
\bid{issn={0090-5364}, mr={0468014}}
\bptok{imsref}%
\end{barticle}
\endbibitem

\bibitem[\protect\citeauthoryear{Stanford and
  Raftery}{2002}]{Stanford-Raftery02}
\begin{barticle}[author]
\bauthor{\bsnm{Stanford},~\bfnm{D.~C.}\binits{D.~C.}} \AND
  \bauthor{\bsnm{Raftery},~\bfnm{A.~E.}\binits{A.~E.}}
(\byear{2002}).
\btitle{Approximate Bayes factors for image segmentation: The pseudolikelihood
  information criterion (PLIC)}.
\bjournal{IEEE Transactions on Pattern Analysis and Machine Intelligence}
\bvolume{24}
\bpages{1517--1520}.
\bptok{imsref}%
\end{barticle}
\endbibitem

\bibitem[\protect\citeauthoryear{Wang, Ju and Wang}{2009}]{Wang-et-al09}
\begin{barticle}[mr]
\bauthor{\bsnm{Wang},~\bfnm{Jie}\binits{J.}},
  \bauthor{\bsnm{Ju},~\bfnm{Lili}\binits{L.}} \AND
  \bauthor{\bsnm{Wang},~\bfnm{Xiaoqiang}\binits{X.}}
(\byear{2009}).
\btitle{An edge-weighted centroidal {V}oronoi tessellation model for image
  segmentation}.
\bjournal{IEEE Trans. Image Process.}
\bvolume{18}
\bpages{1844--1858}.
\bid{doi={10.1109/TIP.2009.2021087}, issn={1057-7149}, mr={2750696}}
\bptok{imsref}%
\end{barticle}
\endbibitem

\bibitem[\protect\citeauthoryear{Yao}{1988}]{Yao88}
\begin{barticle}[mr]
\bauthor{\bsnm{Yao},~\bfnm{Yi-Ching}\binits{Y.-C.}}
(\byear{1988}).
\btitle{Estimating the number of change-points via {S}chwarz' criterion}.
\bjournal{Statist. Probab. Lett.}
\bvolume{6}
\bpages{181--189}.
\bid{doi={10.1016/0167-7152(88)90118-6}, issn={0167-7152}, mr={0919373}}
\bptok{imsref}%
\end{barticle}
\endbibitem

\bibitem[\protect\citeauthoryear{Zhang and Modestino}{1990}]{Zhang-Modestino90}
\begin{barticle}[author]
\bauthor{\bsnm{Zhang},~\bfnm{J.}\binits{J.}} \AND
  \bauthor{\bsnm{Modestino},~\bfnm{J.~W.}\binits{J.~W.}}
(\byear{1990}).
\btitle{A model-fitting approach to cluster validation with application to
  stochastic model-based image segmentation}.
\bjournal{IEEE Transactions on Pattern Analysis and Machine Intelligence}
\bvolume{12}
\bpages{1009--1017}.
\bptok{imsref}%
\end{barticle}
\endbibitem

\bibitem[\protect\citeauthoryear{Zhu and Yuille}{1996}]{Zhu-Yuille96}
\begin{barticle}[author]
\bauthor{\bsnm{Zhu},~\bfnm{Song~Chun}\binits{S.~C.}} \AND
  \bauthor{\bsnm{Yuille},~\bfnm{Alan}\binits{A.}}
(\byear{1996}).
\btitle{Region competition: Unifying snakes, region growing, and {B}ayes/{MDL}
  for multiband image segmentation}.
\bjournal{IEEE Transactions on Pattern Analysis and Machine Intelligence}
\bvolume{18}
\bpages{884--900}.
\bptok{imsref}%
\end{barticle}
\endbibitem

\end{thebibliography}
\end{document}